\documentclass[journal,twoside,web]{ieeecolor}

\def\BibTeX{{\rm B\kern-.05em{\sc i\kern-.025em b}\kern-.08em T\kern-.1667em\lower.7ex\hbox{E}\kern-.125emX}} 

\usepackage{generic}
\usepackage{cite}
\usepackage{amsmath,amssymb,amsfonts}
\usepackage{graphicx}

\usepackage{textcomp}

\usepackage{mathtools}
\usepackage{color}

\usepackage{tikz}
\usetikzlibrary{shapes,arrows}
\usepackage{pgfplots} 
\usepackage{pgfgantt}
\usepackage{pdflscape}
\usepackage{nicefrac}
\usepackage{prettyref}
\pgfplotsset{compat=1.5}
\pgfplotsset{plot coordinates/math parser=false}
\newlength\fwidth
\usepackage{grffile}
\usetikzlibrary{plotmarks}
\usetikzlibrary{arrows.meta}
\usepgfplotslibrary{patchplots}
\usepackage{units}
\usepackage{enumerate}
\let\labelindent\relax
\usepackage{enumitem}
\usepackage{algorithm}
\usepackage{algpseudocode}

\newtheorem{definition}{Definition}
\newtheorem{remark}{Remark}

\newtheorem{theorem}{Theorem}

\newtheorem{problem}{Problem}
\newtheorem{proposition}{Proposition}
\newtheorem{assumption}{Assumption}

\newcommand{\defeq}{\vcentcolon=}
\newcommand{\T}{\scriptscriptstyle\top}       
\newcommand{\mathmin}{\operatorname*{min}}

\newcommand{\mathst}{\text{s.t.}}

\DeclareMathOperator*{\argmax}{arg\,\operatorname*{max}}
\DeclareMathOperator*{\argmin}{arg\,\operatorname*{min}}

\definecolor{myred}{RGB}{233,72,73}%
\definecolor{mygreen}{RGB}{113,191,110}%
\definecolor{myblue}{RGB}{93,147,191}%
\definecolor{mydarkblue}{RGB}{57,101,181}%
\definecolor{mycyan}{rgb}{ 0.05, 0.80, 0.75}%
\definecolor{mygray}{RGB}{163,163,163}%
\definecolor{myorange}{RGB}{233,180,73}%

\long\def\comment#1{}

\algnewcommand{\algorithmicgoto}{\textbf{go to}}%
\algnewcommand{\Goto}[1]{\algorithmicgoto~\ref{#1}}
\usepackage{caption}
\usepackage{subcaption}

\ifdefined\DeclareUnicodeCharacter
  \DeclareUnicodeCharacter{2212}{-}
\fi

\newcommand*\mystrut[1]{\vrule width0pt height0pt depth#1\relax}
\makeatletter
\def\widebreve{\mathpalette\wide@breve}
\def\wide@breve#1#2{\sbox\z@{$#1#2$}%
     \mathop{\vbox{\m@th\ialign{##\crcr
\kern0.18em\vspace{-0.016cm}\brevefill#1{0.55\wd\z@}\crcr\noalign{\nointerlineskip}%
                    $\hss#1#2\hss$\crcr}}}\limits}
\def\brevefill#1#2{$\m@th\sbox\tw@{$#1($}%
  \hss\resizebox{#2}{\wd\tw@}{\rotatebox[origin=c]{90}{\upshape(}}\hss$}
\makeatother

\ifdefined\pdfminorversion\pdfminorversion=4\fi

\makeatletter
\renewcommand\paragraph{%
  \@startsection{paragraph}{4}{\z@}%
    {1ex \@plus .2ex}%
    {-0.3em}
    {\normalfont\normalsize}}
\makeatother

\usepackage{hyperref}
\hypersetup{hidelinks}

\usepackage{eso-pic}
\AddToShipoutPictureBG*{%
\AtPageUpperLeft{%
\setlength\unitlength{1in}%
\hspace*{\dimexpr0.5\paperwidth\relax}
\makebox(0,-0.6)[c]{
\begin{tabular}{c c}

To appear in IEEE Transactions on Automatic Control

\end{tabular}}}}

\AddToShipoutPictureBG*{%
\AtPageUpperLeft{%
\setlength\unitlength{1in}%
\hspace*{\dimexpr0.5\paperwidth\relax}
\makebox(0,-21.35)[c]{
\footnotesize
\begin{tabular}{c c}
© 2025 IEEE. Personal use of this material is permitted. Permission from IEEE must be obtained for all other uses, in any \\
current or future media, including reprinting/republishing this material for advertising or promotional purposes, creating new \\
collective works, for resale or redistribution to servers or lists, or reuse of any copyrighted component of this work in other works.
\end{tabular}}}}

\begin{document}
\title{Verification and Synthesis of Discrete-Time Control Barrier Functions}
\author{Erfan Shakhesi, \IEEEmembership{Student Member, IEEE}, W. P. M. H. (Maurice) Heemels, \IEEEmembership{Fellow, IEEE}, \mbox{and Alexander Katriniok}, \IEEEmembership{Senior Member, IEEE}
\thanks{The research leading to these results has received funding from the European Research Council under the Advanced ERC Grant Agreement PROACTHIS, no. 101055384.}
\thanks{The authors are with the Control Systems Technology Section, Mechanical Engineering, Eindhoven University of Technology, 5600 MB Eindhoven, The Netherlands (e-mail: \href{mailto:e.shakhesi@tue.nl}{e.shakhesi@tue.nl}; 
\href{mailto:m.heemels@tue.nl}{m.heemels@tue.nl}; 
\href{mailto:a.katriniok@tue.nl}{a.katriniok@tue.nl}).
}}

\maketitle

\begin{abstract}
Discrete-time Control Barrier Functions (DTCBFs) have recently attracted interest for guaranteeing safety and synthesizing safe controllers for discrete-time dynamical systems. This paper addresses the open challenges of verifying candidate DTCBFs and synthesizing DTCBFs for general nonlinear discrete-time systems with input constraints and arbitrary safe sets. In particular, we propose a branch-and-bound method, inspired by the $\alpha$BB algorithm, for the verification of candidate DTCBFs in both cases, whether a corresponding control policy is known or unknown. We prove that this method, in a finite number of iterations, either verifies a given candidate function as a valid DTCBF or falsifies it by providing a counterexample (within predefined tolerances). As a second main contribution, we propose a novel bilevel optimization approach to synthesize a DTCBF and a corresponding control policy in finite time. This involves determining the unknown coefficients of a parameterized DTCBF and a parameterized control policy. Furthermore, we introduce various strategies to reduce the computational burden of the bilevel approach. We also demonstrate our methods using numerical case studies.
\end{abstract}


\section{Introduction} \label{sec:introduction}
\label{sec:sec1}
\IEEEPARstart{S}{afety}-critical systems, such as those in aerospace, healthcare, and automotive applications, require formal safety guarantees to prevent operational failures that may lead to substantial harm or damage \cite{Knight2002}. Therefore, the design of controllers that guarantee safety for such systems is crucial.

In the context of dynamical systems, safety is often defined by confining system trajectories within a predefined safe set. Under this definition, safety can be guaranteed by constructing a controlled invariant set, preferably sufficiently large, within the safe set. To achieve this and to synthesize safe controllers, one notable approach discussed in the literature is the use of Control Barrier Functions (CBFs) \cite{Ames2014a}. In particular, CBF-based controllers often act as safety filters, adjusting nominal control inputs when these inputs would lead to a violation of safety constraints. This control technique was initially applied to continuous-time systems, with application examples including adaptive cruise control \cite{Ames2014a}, bipedal robots \cite{Nguyen2020}, and multi-agent systems \cite{Wang2017a}. An extension of CBFs to discrete-time dynamical systems, referred to as discrete-time CBFs (DTCBFs), has been introduced in \cite{Agrawal2017a}. DTCBFs are also used within Model Predictive Control (MPC) both to formulate stage constraints \cite{Zeng2021a, Zeng2021b} and to guarantee recursive feasibility \cite{katriniok2023}, and hence have broad application potential.

Synthesizing (DT)CBFs for general systems subject to input constraints is inherently challenging and is currently an open research question in the literature. In this paper, we aim to address this challenge for the discrete-time case.
\subsection{Related Work}
As (DT)CBFs can be regarded as a generalization of \textit{classical} set invariance, we review previous work on synthesizing controlled invariant sets, as well as on verifying candidate (DT)CBFs and synthesizing (DT)CBFs. Note that the term ``candidate'' refers to any function without specific requirements, which may or may not satisfy the (DT)CBF conditions and thus may or may not be a valid (DT)CBF.
Our primary focus is on discrete-time systems, but we also discuss existing work on continuous-time systems, as some of the ideas could be useful in discrete-time cases as well.
\subsubsection{Continuous-Time Systems -- Control Barrier Functions}
    Using CBFs to guarantee safety offers several advantages, including the simplicity of synthesizing safe controllers and the smooth adjustment of nominal control inputs \cite{Wabersich2023}. A major drawback is the challenge of synthesizing CBFs, as already mentioned. Common methods for the synthesis of CBFs are Sum-of-Squares (SOS) programming \cite{Ames2019a, Wang2023a, clark2022, dai2022, clark2021, Kang2023}, backup policy methods \cite{chen2021, Gurriet2020}, and learning-based methods \cite{wei2023, Zhang2023, lindemann2024, tan2022, verification_neural_CBF}:
    
        \textit{SOS Programming:} The conditions required to verify or synthesize polynomial CBFs for polynomial systems and semi-algebraic safe sets can be mapped to SOS constraints via the Positivstellensatz \cite{Ames2019a}. However, due to the bilinearity introduced by the multiplication of unknown polynomials, these constraints cannot be imposed directly in SOS programming. To address this issue and enable the synthesis of CBFs for control-affine polynomial systems, alternating-descent approaches based on SOS programming are proposed in \cite{Wang2023a, clark2022, dai2022, Ames2019a}. Regarding the verification of polynomial candidate CBFs for control-affine polynomial systems, SOS programming is used in \cite{dai2022} and \cite{clark2021}. Additionally, \cite{Kang2023} uses optimization problems with SOS programming to verify and synthesize \textit{robust} CBFs for control-affine polynomial systems with convergence guarantees.

        It is noteworthy that these methods cannot be used for the verification or synthesis of DTCBFs for discrete-time systems. This limitation arises because the DTCBF constraint is not affine in the control input. As a result, the bilinear terms in the DTCBF constraint are distinct from, and more complex than, those in the continuous-time case, and cannot be addressed by the proposed alternating-descent approaches \cite{Shakhesi2025-SOS}.
        
        \textit{Backup Policy Methods:} 
        In \cite{chen2021} and \cite{Gurriet2020}, CBFs are synthesized by enlarging an initial controlled invariant set through forward integration of dynamics under a known stabilizing controller, referred to as the backup policy. The challenge with this method lies in the complexity of predetermining both a stabilizing control policy and an initial controlled invariant set for an arbitrary system.
        
        \textit{Learning-Based Methods:}
        Candidate CBFs can be synthesized by generating functions that adhere to the CBF conditions at a finite set of sample states using an evolutionary algorithm \cite{wei2023}, mixed-integer convex programming with SOS programming for polynomial CBFs considering polynomial systems \cite{Zhang2023}, or convex optimization programming \cite{lindemann2024}. Under additional conditions such as the density of the sample states and the smoothness of the system dynamics, these candidate CBFs are guaranteed to be valid CBFs \cite{wei2023, Zhang2023, tan2022, lindemann2024}. Additionally, regarding the verification of candidate CBFs, \cite{tan2022} proposes an optimization-based method with grid sampling, and \cite{verification_neural_CBF} introduces a framework for ReLU-based neural CBFs through symbolic derivative bound propagation.         \vspace{1em}
        
        \subsubsection{Discrete-Time Systems}
        \paragraph{Classical Numerical Methods for Controlled Invariant Sets}
        
        \hspace*{0.7em}\textit{Linear Systems:} Linear programming is used in \cite{blanchini2008} to synthesize polyhedral and ellipsoidal controlled invariant sets for linear systems using linear controllers.
        
        \textit{Polynomial Systems:} \cite{Wang2023a} presents an iterative method for control-affine polynomial systems based on SOS programming that starts from a relatively small predetermined controlled invariant set to synthesize polyhedral controlled invariant sets. This method, however, does not guarantee convergence, and the properties (such as shape and size) of the resulting set depend strongly on the predetermined controlled invariant set. Additionally, \cite{Korda2014} proposes an SOS programming-based method to compute an outer approximation of the maximal controlled invariant (MCI) set for polynomial systems. Such an outer approximation, however, is not suitable for guaranteeing safety, since it is not necessarily controlled invariant and may include unsafe states.
        
        \textit{General Systems:} For constrained nonlinear control-affine systems, \cite{CANNON20031487} uses linear programming to synthesize low-complexity polytopic controlled invariant sets with \textit{linear controllers}, approximating nonlinear systems with Linear Differential Inclusions (LDIs). However, the resulting controlled invariant sets may be conservative for the original nonlinear systems. Moreover, \cite{Wan2009} exploits a numerical approach via interval analysis to synthesize relatively more complex polyhedral controlled invariant sets using \textit{linear controllers} assumed to be \textit{a priori} given. Although these methods are of interest, the main drawback is that assigning linear controllers is not always a valid choice for synthesizing polytopic controlled invariant sets, even for simple systems \cite{Alessandro2001}.
   
    \paragraph{Reachability-Based Methods}
    One common approach to synthesizing controlled invariant sets for discrete-time dynamical systems is to compute backward reachable sets. By using such a set in a recursive way, MCI sets can be obtained upon termination \cite{BERTSEKAS1972,blanchini2008, Gilbert1991, Rungger2017}. However, finite termination is not guaranteed for general systems. In a finite number of iterations, if the procedure starts from a controlled invariant set, an inner approximation of the MCI set is computed, and if it starts from the state constraint set, an outer approximation of the MCI set is obtained \cite{Rungger2017}. To compute one-step backward reachable sets, \cite{Gilbert1991} uses optimization-based methods for linear systems, and branch-and-bound approaches are proposed in \cite{BRAVO2005} and \cite{brown2023} for nonlinear systems. It should be noted that these methods can be computationally expensive, and the resulting sets may not have closed-form expressions, posing a challenge for controller synthesis (e.g., implementation in MPC frameworks).
    
    In addition to these works, \cite{scalablesynthesisformallyverified} proposes a learning-based method to synthesize verified neural value functions for Hamilton-Jacobi reachability analysis, assuming a given controller. While the method is interesting and computationally scalable, the assumption that the control policy is known is restrictive.
    \paragraph{Discrete-Time Control Barrier Functions}
    To the best of our knowledge, there is no systematic method available in the literature for the direct synthesis of DTCBFs (specifically with termination guarantees). Notably, prior works that exploited DTCBFs, such as \cite{Agrawal2017a}, \cite{Zeng2021a}, and \cite{Zeng2021b}, have resorted to using \textit{candidate} DTCBFs. These candidates, however, may not rigorously adhere to the formal DTCBF definition, particularly in the presence of input constraints, thereby not providing direct safety and feasibility guarantees. Nevertheless, \cite{Zhang2024} proposes a learning-based method for synthesizing candidate DTCBFs by minimizing safety violations on sample states. Additionally, \cite{Freire2023a} uses the maximal output admissible set (MOAS) theory for synthesizing a particular type of DTCBF (considering the class-$\mathcal{K}_{\infty}$ function as the identity function in the formal DTCBF definition). Specifically, the MOAS is the set of all initial states for which the resulting trajectories, under an \textit{a priori} given stabilizing control policy, satisfy state and input constraints. The MOAS theory, though, has not been well-studied for nonlinear systems. Moreover, \cite{Shakhesi2025} proposes a counterexample-guided approach for synthesizing DTCBFs. In this approach, a candidate DTCBF and a corresponding control policy are trained to satisfy the required conditions on a set of sample states. A verification algorithm then attempts to find a counterexample. If one is found, it is added to the set of sample states to update the candidate DTCBF and the control policy. This procedure continues until a valid DTCBF is obtained. While effective in many practical examples, this approach does not guarantee finding a valid DTCBF in finite time or ensure that its zero-superlevel set is maximal in size.
\subsection{Main Contributions}
Based on the above literature review, synthesizing controlled invariant sets for general nonlinear discrete-time systems, particularly in the case where a corresponding control policy is unknown, and verifying and synthesizing DTCBFs, even for linear systems and especially with termination guarantees, are important open research questions. This paper aims to fill this gap by providing the following key contributions:
\begin{enumerate}[label=(\roman*)]
    \item Regarding the verification of candidate DTCBFs (Section~\ref{sec:sec3}), we present a branch-and-bound (BB) method, and we formally prove that it terminates in a finite number of iterations and either verifies a candidate function as a valid DTCBF or falsifies it by providing a counterexample (within predefined tolerances), in both the known and unknown control policy cases. We already presented this method in our preliminary conference paper \cite{Shakhesi2024}, but without any proofs due to space limitations. The full proofs and additional details are provided in the present paper.
    \item Regarding the synthesis of DTCBFs (Section \ref{sec:sec4}):
    \begin{itemize}
        \item We propose a novel bilevel optimization approach to synthesize a DTCBF and a corresponding control policy for a general nonlinear discrete-time system with input constraints and a safe set. This involves determining the unknown coefficients of a parameterized DTCBF and a parameterized control policy.
        \item We guarantee convergence to the DTCBF with the largest zero-superlevel set (in an appropriate sense) for a discrete-time system and a safe set, provided that there exist coefficients for the parameterized DTCBF and the parameterized control policy that make the DTCBF valid.
        \item We propose various strategies to reduce the computation time of the bilevel approach.
    \end{itemize}
    \item We apply our proposed methods to numerical case studies (Section \ref{sec:sec5}).
\end{enumerate}
\subsection{Notation and Definitions}
We use $\mathbb{R}$, $\mathbb{R}_{>0}$, and $\mathbb{R}_{\geqslant 0}$ to denote the sets of real numbers, positive real numbers, and non-negative real numbers, respectively. Additionally, $\mathbb{R}^{n}$ represents the set of all $n$-dimensional vectors of real numbers. For a vector $x \in \mathbb{R}^n$, $x_i \in \mathbb{R}$, $i \in \{1,2, \hdots, n\}$, represents the \mbox{$i$-th} element of $x$, and $x^{\T}$ represents its transpose. Moreover, $\mathbb{N} \defeq \{1,2,3, \hdots\}$, and $\mathbb{N}_0 \defeq \mathbb{N} \cup \{0\}$. For a set $\mathcal{C} \subset \mathbb{R}^n$, $\partial \mathcal{C}$ represents its boundary.
\begin{definition}[Zero-Superlevel Set]
The \textit{zero-superlevel set} $\mathcal{C}$ of a function $h:\mathbb{R}^{n} \rightarrow \mathbb{R}$ is defined as 
\begin{equation}
    \mathcal{C} \defeq \{ x \in \mathbb{R}^{n} \mid h(x) \geqslant 0 \}. \nonumber
\end{equation}
\end{definition}
\vspace{0.1cm}
\begin{definition}[Class-$\mathcal{K}_\infty$ Functions]
    A continuous function $\gamma:\mathbb{R}_{\geqslant 0} \rightarrow \mathbb{R}_{\geqslant 0}$ is said to be a \textit{class-$\mathcal{K}_\infty$ function}, denoted by $\gamma \in \mathcal{K}_{\infty}$, if it is strictly increasing, $\gamma(0) = 0$, and $\lim_{r\rightarrow \infty} \gamma(r)=\infty$.
\end{definition}

For $\gamma \in \mathcal{K}_{\infty}$, we use the notation $\gamma \in \mathcal{K}^{\leqslant \mathrm{id}}_{\infty}$ to indicate that $\gamma(r) \leqslant r$ for all $r \in \mathbb{R}_{\geqslant 0}$.
\begin{definition}[$n$-Rectangle Set]
   We define $X \subset \mathbb{R}^n$ as an \textit{$n$-rectangle set}, using the notation $X \defeq [x^{lb}, ~ x^{ub}]$, where $x^{lb},x^{ub} \in \mathbb{R}^n$ with $x_i^{lb} \leqslant x_i^{ub}$, $i \in \{1,2, \hdots, n\}$, given by 
    \begin{align}
        X \defeq \{ x \in \mathbb{R}^n &\mid x_i^{lb} \leqslant x_i \leqslant x_i^{ub}, ~i \in \{1,2, \hdots, n\}\}. \nonumber
    \end{align}
\end{definition}

\section{Background on Control Barrier Functions}
\label{sec:sec2}
Consider the discrete-time system
\begin{equation} \label{eq:sec2:dynamical-system}
    x_{t+1} = f(x_t, u_t),
\end{equation}
with state vector $x_t \in \mathbb{R}^{n}$, control input vector \mbox{$u_t \in \mathbb{U} \subseteq \mathbb{R}^{m}$}, both at discrete time $t \in \mathbb{N}_0$, and mapping \mbox{$f:\mathbb{R}^{n} \times \mathbb{R}^{m} \rightarrow \mathbb{R}^{n}$}. Here, $\mathbb{U}$ is the control admissible set. For the system \eqref{eq:sec2:dynamical-system}, the safe set $\mathcal{S}$ is defined as
\begin{align} \label{eq:sec2:safe-set}
    \mathcal{S} \defeq \{ x \in \mathbb{R}^{n} \mid s(x) \geqslant 0 \},
\end{align}
where $s:\mathbb{R}^{n} \rightarrow \mathbb{R}$ is a given mapping.
\begin{definition}[Controlled Invariance \cite{Ames2019a}]
For the system \eqref{eq:sec2:dynamical-system} with the control admissible set $\mathbb{U}$, a set $\mathcal{C} \subset \mathbb{R}^{n}$ is \textit{controlled invariant} if, for each $x \in \mathcal{C}$, there exists a control input $u \in \mathbb{U}$ such that $f(x,u) \in \mathcal{C}$.
\end{definition}
\begin{definition}[Safety] \label{def:sec2:safety} 
    The system \eqref{eq:sec2:dynamical-system} with the control admissible set $\mathbb{U}$ is considered safe with respect to a safe set \mbox{$\mathcal{S} \subseteq \mathbb{R}^n$} and an initial state set $X_0 \subseteq \mathcal{S}$, denoted by $(\mathcal{S}, X_0)$-\textit{safe}, if there exists a control policy $\mu:\mathcal{S} \rightarrow \mathbb{U}$ such that $x_t \in \mathcal{S}$ for all $t \in \mathbb{N}_0$ and all trajectories of $x_{t+1} = f(x_t, \mu(x_t)) \in \mathcal{S}$ with $x_0 \in X_0$.
\end{definition}

Based on Definition \ref{def:sec2:safety}, one way to guarantee \textit{$(\mathcal{S}, X_0)$-safety} for the system \eqref{eq:sec2:dynamical-system} with the control admissible set $\mathbb{{U}}$ is through constructing a \textit{controlled invariant} set $\mathcal{C}$ such that \mbox{$X_0 \subseteq \mathcal{C} \subseteq \mathcal{S}$}. To achieve this, a method that has gained increasing attention in the literature is the use of DTCBFs \cite{Agrawal2017a}, \cite{Zeng2021a}, which also enables the synthesis of safe controllers.
\begin{definition}[DTCBF \cite{Agrawal2017a}, \cite{Zeng2021a}] \label{def:sec2:DTCBF}
     Consider a function \mbox{$h:\mathbb{R}^n \rightarrow \mathbb{R}$} with zero-superlevel set $\mathcal{C}$. For the system~\eqref{eq:sec2:dynamical-system} with the control admissible set $\mathbb{U}$, $h$ is a \textit{discrete-time Control Barrier Function (DTCBF)} if there exists a $\gamma \in \mathcal{K}^{\leqslant \mathrm{id}}_{\infty}$ such that, for each $x \in \mathcal{C}$, there is a control input $u \in \mathbb{U}$ satisfying
    \begin{align} \label{eq:sec2:DTCBF-constraint-alpha}
        h(f(x,u)) - h(x) \geqslant - \gamma(h(x)).
    \end{align}
    We refer to \eqref{eq:sec2:DTCBF-constraint-alpha} as the DTCBF constraint.
\end{definition}
\begin{theorem} \label{theorem:sec2:safety}
    Consider a function $h:\mathbb{R}^n \rightarrow \mathbb{R}$ with zero-superlevel set $\mathcal{C}$. The set $\mathcal{C}$ is controlled invariant for the system~\eqref{eq:sec2:dynamical-system} with the control admissible set $\mathbb{U}$ if and only if $h$ is a DTCBF. Moreover, the system~\eqref{eq:sec2:dynamical-system} with $\mathbb{U}$ is \textit{$(\mathcal{S}, X_0)$-safe} if $X_0 \subseteq \mathcal{C} \subseteq \mathcal{S}$.
    
\proof
    Suppose $h$ is a DTCBF. Then there exist a \mbox{$\gamma \in \mathcal{K}^{\leqslant \mathrm{id}}_{\infty}$} and a control policy $\pi:\mathcal{C} \rightarrow \mathbb{U}$ such that, for all $x \in \mathcal{C}$, 
    \begin{align}
        h(f(x,\pi(x))) \geqslant h(x) - \gamma(h(x)) \geqslant 0. \nonumber
    \end{align}
    Hence, if $x \in \mathcal{C}$, there exists an input $u := \pi(x) \in \mathbb{U}$ such that $f(x, u) \in \mathcal{C}$, and thus $\mathcal{C}$ is controlled invariant. To prove the converse, assume that $\mathcal{C}$ is controlled invariant. Then there exists a control policy $\pi: \mathcal{C} \rightarrow \mathbb{U}$ such that for all $x \in \mathcal{C}$, 
    \begin{align}
        h(f(x,\pi(x))) \geqslant 0. \nonumber
    \end{align}
    Thus, the DTCBF constraint \eqref{eq:sec2:DTCBF-constraint-alpha} holds with $\gamma(r) = r$, \mbox{$r \in \mathbb{R}_{\geqslant 0}$}, which belongs to $\mathcal{K}_{\infty}^{\leqslant \mathrm{id}}$. Therefore, $h$ is a DTCBF. Moreover, according to Definition \ref{def:sec2:safety}, if $X_0 \subseteq \mathcal{C} \subseteq \mathcal{S}$, the system \eqref{eq:sec2:dynamical-system} with $\mathbb{U}$ is \textit{$(\mathcal{S}, X_0)$-safe}.
\endproof
\end{theorem}

In addition to the fact that the zero-superlevel set $\mathcal{C}$ of a DTCBF $h$ is \textit{controlled invariant} (see Theorem \ref{theorem:sec2:safety}), the DTCBF constraint \eqref{eq:sec2:DTCBF-constraint-alpha} incorporates an additional term, represented by $\gamma \in \mathcal{K}^{\leqslant \mathrm{id}}_{\infty}$. This term regulates the rate at which the states of the system \eqref{eq:sec2:dynamical-system} can approach the boundary of $\mathcal{C}$. The advantages of adjusting $\gamma$ are discussed in more detail in \cite{Zeng2021a} and \cite{katriniok2023}.
 
For simplicity in this paper, we assume that $\gamma \in \mathcal{K}^{\leqslant \mathrm{id}}_{\infty}$ in Definition~\ref{def:sec2:DTCBF} is \textit{a priori} given and fixed, similar to \cite{chen2021, tonkens2022refining, tan2022}, and, in the remainder of this paper, we say that $h$ with a given $\gamma$, denoted by $(h,\gamma)$, is a DTCBF for the system \eqref{eq:sec2:dynamical-system} with $\mathbb{U}$ if, for all $x \in \mathcal{C}$, there exists an input $u \in \mathbb{U}$ such that the DTCBF constraint \eqref{eq:sec2:DTCBF-constraint-alpha} is satisfied. Moreover, we say that a control policy \mbox{$\pi:\mathcal{C} \rightarrow \mathbb{U}$} is \textit{a friend of} $(h, \gamma)$, assuming $\mathbb{U}$ and $f$ are clear from the context, if, for all $x \in \mathcal{C}$,
\begin{align} 
    \pi(x) \in \mathbb{U} ~~~ \text{and} ~~~ h(f(x,\pi(x))) - h(x) \geqslant - \gamma(h(x)). \nonumber
\end{align}
However, it is worth noting that the methods proposed in this paper remain applicable even in cases where $\gamma$ is unknown (see Remark~\ref{remark:sec3:unknown-gamma} for the verification problem and Remark~\ref{sec4:remark-existence-class-kappa} for the synthesis problem below). 

\section{Verification}
\label{sec:sec3}
In this section, we will provide methods for the verification of candidate DTCBFs, and in Section \ref{sec:sec4}, we will address the problem of synthesizing DTCBFs.
\subsection{Problem Statement and Assumptions} \label{sec:sec3:problem-statement}
\begin{problem}[Verification]
\label{problem:verification}
   Consider a given candidate DTCBF \mbox{$h:\mathbb{R}^n \rightarrow \mathbb{R}$} and a given $\gamma \in \mathcal{K}^{\leqslant \mathrm{id}}_{\infty}$. The objective is to either verify $(h, \gamma)$ as a valid DTCBF for the system \eqref{eq:sec2:dynamical-system} with the control admissible set $\mathbb{U}$, or falsify it by providing a counterexample, in case a control policy $\pi:\mathcal{C} \rightarrow \mathbb{U}$ is 
   \begin{enumerate}[label=(\roman*), ref=Case (\roman*)]
       \item \textit{a priori} given, or \label{case:problem-verification:known}
       \item unknown. \label{case:problem-verification:unknown}
   \end{enumerate}
\end{problem}
\begin{assumption} \label{assump:sec3:C-bounded}
    The zero-superlevel set $\mathcal{C}$ of the candidate DTCBF $h$ and the control admissible set $\mathbb{U}$ are compact. 
\end{assumption}
\begin{assumption} \label{assump:sec3:continuous-functions}
    The mapping $f$ associated with the system~\eqref{eq:sec2:dynamical-system}, the candidate DTCBF $h$, and the control policy $\pi$, if given, are continuous.
\end{assumption}
\begin{remark} \label{remark:sec3:unknown-gamma}
    If $\gamma \in \mathcal{K}^{\leqslant \mathrm{id}}_{\infty}$ is not \textit{a priori} given, we can still apply our BB method. Specifically, by considering $\gamma$ as the identity function, the DTCBF constraint \eqref{eq:sec2:DTCBF-constraint-alpha} represents the least restrictive condition compared to other choices of $\gamma \in \mathcal{K}^{\leqslant \mathrm{id}}_{\infty}$. Hence, we can assert that $h$ is a DTCBF if and only if $h$, with $\gamma$ defined as the identity function, is a valid DTCBF. Then, a desired $\gamma$ for a valid DTCBF can be determined using, e.g., systematic trial-and-error approaches.
\end{remark}
\subsection{Optimization Framework} \label{sec:sec3:optimization-framework}
\begin{proposition}[Known Control Policy] \label{prop:sec3:verification-known}
    Consider a candidate DTCBF $h:\mathbb{R}^n \rightarrow \mathbb{R}$, its zero-superlevel set $\mathcal{C}$, a $\gamma \in \mathcal{K}^{\leqslant \mathrm{id}}_{\infty}$, and a control policy $\pi:\mathcal{C} \rightarrow \mathbb{U}$. Additionally, consider an $n$-rectangle set $\mathbb{X}$ such that $\mathcal{C} \subseteq \mathbb{X}$. For the system~\eqref{eq:sec2:dynamical-system} with the control admissible set $\mathbb{U}$, $(h,\gamma)$ is a DTCBF and $\pi$ is a friend of $(h, \gamma)$ if and only if $\mathcal{F}^* \geqslant 0$, where $\mathcal{F}^* \in \mathbb{R}$ is the global minimum of
    \begin{subequations} \label{eq:sec3:known-proposition}
    \begin{align}
        \mathcal{F}^* \defeq \min_{x \in \mathbb{X}} ~&~ h(f(x, \pi(x))) - h(x) + \gamma(h(x)) \label{eq:sec3:known-proposition-objective} \\
        \mathst ~&~ -h(x) \leqslant 0. \label{eq:sec3:known-proposition-constraint}
    \end{align}
    \end{subequations}
    \begin{proof}
    By Assumption~\ref{assump:sec3:C-bounded}, the set $\mathcal{C}$ is compact, and therefore, there exists an $n$-rectangle set $\mathbb{X}$ such that \mbox{$\mathcal{C} \subseteq \mathbb{X}$}.
    Suppose first that $\mathcal{F}^* \geqslant 0$. Then, for all $x \in \mathcal{C}$, the DTCBF constraint \eqref{eq:sec2:DTCBF-constraint-alpha} is satisfied. Hence, by Definition~\ref{def:sec2:DTCBF}, $(h, \gamma)$ is a DTCBF, and $\pi$ is a friend of $(h, \gamma)$ for the system~\eqref{eq:sec2:dynamical-system} with $\mathbb{U}$. To prove the converse, we first show that the global minimum $\mathcal{F}^*$ of \eqref{eq:sec3:known-proposition} exists. By the Weierstrass theorem~\cite{ARORA2017}, since the feasible set of the optimization problem \eqref{eq:sec3:known-proposition} is compact (by Assumption~\ref{assump:sec3:C-bounded}) and the objective function \eqref{eq:sec3:known-proposition-objective} is continuous (as $h$, $f$, and $\pi$ are continuous by Assumption~\ref{assump:sec3:continuous-functions} and \mbox{$\gamma \in \mathcal{K}^{\leqslant \mathrm{id}}_{\infty}$}), $\mathcal{F}^*$ exists. Hence, if $(h, \gamma)$ is a DTCBF and $\pi$ is a friend of $(h, \gamma)$, then the DTCBF constraint \eqref{eq:sec2:DTCBF-constraint-alpha} is satisfied for all $x \in \mathcal{C}$, which implies that $\mathcal{F}^* \geqslant 0$.
    \end{proof}
\end{proposition}
\begin{proposition}[Unknown Control Policy] \label{prop:sec3:verification-unknown}
    Consider a candidate DTCBF $h:\mathbb{R}^n \rightarrow \mathbb{R}$, its zero-superlevel set $\mathcal{C}$, and a \mbox{$\gamma \in \mathcal{K}^{\leqslant \mathrm{id}}_{\infty}$}. Additionally, consider an $n$-rectangle set $\mathbb{X}$ such that $\mathcal{C} \subseteq \mathbb{X}$. For the system \eqref{eq:sec2:dynamical-system} with $\mathbb{U}$, $(h, \gamma)$ is a DTCBF if and only if $\mathcal{F}^* \geqslant 0$, where $\mathcal{F}^* \in \mathbb{R}$ is the global optimum of
    \begin{subequations} \label{eq:sec3:unknown-proposition}
    \begin{align}
        \mathcal{F}^* \defeq \min_{x \in \mathbb{X}} \max_{u \in \mathbb{U}} ~&~ h(f(x, u)) - h(x) + \gamma(h(x)) \label{eq:sec3:unknown-proposition-objective} \\
        \mathst ~&~ -h(x) \leqslant 0.
    \end{align}
    \end{subequations}
    \proof
    First, by Assumption~\ref{assump:sec3:C-bounded}, the set $\mathcal{C}$ is compact, and therefore, there exists an $n$-rectangle set $\mathbb{X}$ such that \mbox{$\mathcal{C} \subseteq \mathbb{X}$}. The inner maximization problem in \eqref{eq:sec3:unknown-proposition} corresponds to finding a suitable control input $u \in \mathbb{U}$ that satisfies the DTCBF constraint \eqref{eq:sec2:DTCBF-constraint-alpha} for a given $x \in \mathcal{C}$. The outer minimization problem checks this condition for all $x \in \mathcal{C}$. As a result, $\mathcal{F}^* \geqslant 0$ implies that for all $x \in \mathcal{C}$, there exists a control input $u \in \mathbb{U}$ such that the DTCBF constraint \eqref{eq:sec2:DTCBF-constraint-alpha} is satisfied. Thus, $(h,\gamma)$ is a valid DTCBF. To prove the converse, we first show that $\mathcal{F}^*$ exists. The global maximum of the inner problem for a fixed $x \in \mathcal{C}$ exists by the Weierstrass theorem~\cite{ARORA2017}, since $\mathbb{U}$ is compact and the objective function \eqref{eq:sec3:unknown-proposition-objective} is continuous by Assumptions~\ref{assump:sec3:C-bounded}~and~\ref{assump:sec3:continuous-functions}. The global minimum of the outer problem also exists, as $\mathcal{C}$ is compact and the objective function~\eqref{eq:sec3:unknown-proposition-objective} is continuous. Hence, if $(h, \gamma)$ is a DTCBF, then for all $x \in \mathcal{C}$, there exists a control input $u \in \mathbb{U}$ satisfying the DTCBF constraint \eqref{eq:sec2:DTCBF-constraint-alpha}. Therefore, $\mathcal{F}^* \geqslant 0$. 
    \endproof
\end{proposition}
As outlined in Propositions \ref{prop:sec3:verification-known} and \ref{prop:sec3:verification-unknown}, to verify whether $(h, \gamma)$ is a DTCBF, one can verify or falsify the non-negativity of the global optimum $\mathcal{F}^*$ of \eqref{eq:sec3:known-proposition} or \eqref{eq:sec3:unknown-proposition}. To this end and since these optimization problems are generally non-convex, we aim to propose a BB method, inspired by the $\alpha$BB algorithm \cite{Floudas2010a}.

\begin{remark}
    If we also aim to verify $\mathcal{C} \subseteq \mathcal{S}$ (see Theorem~\ref{theorem:sec2:safety}), that is, to check whether the system~\eqref{eq:sec2:dynamical-system} with $\mathbb{U}$ is \textit{$(\mathcal{S}, X_0)$-safe} for some $X_0 \subseteq \mathcal{C}$, we can check whether $s^* \geqslant 0$, where
    \begin{subequations} \label{eq:verification:remark:CsubseteqS}
        \begin{align} 
        s^* \defeq \min_{x \in \mathbb{X}} ~&~ s(x) \\
        \mathst ~&~ -h(x) \leqslant 0.
        \end{align}
    \end{subequations}
    Here, $s:\mathbb{R}^n \rightarrow \mathbb{R}$ is defined in \eqref{eq:sec2:safe-set}. This verification can also be done using the BB method discussed in this paper.
\end{remark}
\subsection{Preliminaries} \label{sec:sec3:preliminaries}
We first recall the definition and the construction of convex underestimators from \cite{Floudas2010a}.
\begin{definition}[Convex Underestimator \cite{Floudas2010a}]
    A function $\widebreve{\mathcal{F}}:\mathbb{R}^n \rightarrow \mathbb{R}$ is called a \textit{convex underestimator} of a given function $\mathcal{F}:\mathbb{R}^n \rightarrow \mathbb{R}$ on a set $X \subseteq \mathbb{R}^n$ if $\widebreve{\mathcal{F}}$ is convex on $X$ and it holds that $\widebreve{\mathcal{F}}(x) \leqslant \mathcal{F}(x)$ for all $x \in X$.
\end{definition}

In the literature, a common approach to construct a convex underestimator of a function $\mathcal{F}:\mathbb{R}^n \rightarrow \mathbb{R}$ on an $n$-rectangle set $X \defeq [x^{lb}, ~ x^{ub}]$ is through adding a quadratic polynomial to $\mathcal{F}$ \cite{Floudas2010a}. This polynomial is chosen such that it is non-positive on $X$ and has sufficiently large coefficients for the quadratic terms to overcome the non-convexitiy of $\mathcal{F}$. To be more precise, the function $\widebreve{\mathcal{F}}:\mathbb{R}^n \rightarrow \mathbb{R}$, constructed as
\begin{align} \label{eq:sec3:preliminaries:convex-underestimator}
    \widebreve{\mathcal{F}}(x) \defeq \mathcal{F}(x) + \sum_{i=1}^{n} \alpha_i\left(x_i^{lb} - x_i\right)\left(x_i^{ub} - x_i\right),
\end{align}
is a convex underestimator of $\mathcal{F}$ on $X \defeq [x^{lb}, ~ x^{ub}]$ if $\alpha_i \in \mathbb{R}_{\geqslant 0}$, $i \in \{1, \hdots, n\}$, are sufficiently large. Numerous techniques for computing suitable values of $\alpha_i$, $i \in \{1, \hdots, n\}$, can be found in \cite{ADJIMAN1998a}, such as the Scaled Gerschgorin method. \\ 
Moreover, two properties of convex underestimators, constructed as in \eqref{eq:sec3:preliminaries:convex-underestimator}, are recalled from \cite{Maranas1994b}.
\begin{proposition}[Maximum Separation \cite{Maranas1994b}] \label{prop:sec3:preliminaries:max-separation}
    The maximum separation between a function $\mathcal{F}:\mathbb{R}^n \rightarrow \mathbb{R}$ and its convex underestimator $\widebreve{\mathcal{F}}:\mathbb{R}^n \rightarrow \mathbb{R}$, as in \eqref{eq:sec3:preliminaries:convex-underestimator}, on an $n$-rectangle set $X \defeq [x^{lb}, \; x^{ub}]$ occurs at the center of $X$ and is given by%
    \begin{align}
        \max_{x\in X} ~ \bigl(\mathcal{F}(x) - \widebreve{\mathcal{F}}(x)\bigr) = \frac{1}{4}\sum_{i=1}^{n} \alpha_i\left(x_i^{ub} - x_i^{lb}\right)^2. \nonumber
    \end{align}
\end{proposition}
\vspace{0.15cm}
\begin{proposition}[Tighter Convex Underestimators \cite{Maranas1994b}] \label{prop:sec3:preliminaries:tighter-underestimator}
Consider two $n$-rectangle sets $X^{(1)} \defeq [x^{lb, (1)}, ~ x^{ub, (1)}]$ and $X^{(2)} \defeq [x^{lb, (2)}, ~ x^{ub, (2)}]$ with $X^{(2)} \subseteq X^{(1)}$. Let $\widebreve{\mathcal{F}}^{(1)}$ and $\widebreve{\mathcal{F}}^{(2)}$ be convex underestimators of a function $\mathcal{F}$ on $X^{(1)}$ and $X^{(2)}$, respectively, constructed as
    \begin{align}
        \widebreve{\mathcal{F}}^{(l)}(x) &\defeq \mathcal{F}(x) + \sum_{i=1}^{n} \alpha_i^{(l)}\bigl(x_i^{lb, (l)} - x_i\bigr)\bigl(x_i^{ub,(l)} - x_i\bigr), \nonumber
    \end{align}
    $l \in \{1,2\}$, where $\alpha_i^{(1)}, \alpha_i^{(2)} \in \mathbb{R}_{\geqslant 0}$ are computed based on one of the methods discussed in \cite{ADJIMAN1998a}.
    Then, it holds that $\alpha_i^{(2)} \leqslant \alpha_i^{(1)}$ for all $i \in \{1,\hdots, n\}$, and $\widebreve{\mathcal{F}}^{(2)}$ is a tighter convex underestimator of $\mathcal{F}$ than $\widebreve{\mathcal{F}}^{(1)}$ on $X^{(2)}$, in the sense that $\widebreve{\mathcal{F}}^{(1)}(x) \leqslant \widebreve{\mathcal{F}}^{(2)}(x) \leqslant \mathcal{F}(x)$ for all $x \in X^{(2)}$.
\end{proposition}
\subsection{Proposed Algorithm for Verification}\label{sec:sec3:verification}
In the case of a known control policy $\pi$ for a candidate DTCBF $h$ with a $\gamma \in \mathcal{K}^{\leqslant \mathrm{id}}_{\infty}$, and given that the optimization problem~\eqref{eq:sec3:known-proposition} is generally non-convex, our algorithm constructs a convex optimization problem using convex underestimators of both the objective function~\eqref{eq:sec3:known-proposition-objective} and the constraint~\eqref{eq:sec3:known-proposition-constraint} on the set $\mathbb{X}$, thereby obtaining a tractable lower bound on the global minimum of \eqref{eq:sec3:known-proposition}. If the global minimum of the convex problem is non-negative, $(h, \gamma)$ with its friend $\pi$ is verified as a valid DTCBF. Otherwise, we partition $\mathbb{X}$ into smaller subdomains to construct tighter convex underestimators. We repeat the same procedure within each subdomain. We approve subdomains in which $(h, \gamma)$ with its friend $\pi$ is verified as a valid DTCBF (i.e., the global minimum of the constructed convex problem is non-negative within the corresponding subdomain) or those that are entirely outside the zero-superlevel set $\mathcal{C}$ of $h$, and we further subdivide the remaining subdomains. We continue this process until either all subdomains are approved, showing that $(h, \gamma)$ with its friend $\pi$ is a valid DTCBF, or a specific $x \in \mathcal{C}$ is found that serves as a counterexample, for which the DTCBF constraint~\eqref{eq:sec2:DTCBF-constraint-alpha} cannot be satisfied by the control policy~$\pi$ (within predefined tolerances).

If a corresponding control policy is unknown, instead of solving the min-max problem \eqref{eq:sec3:unknown-proposition}, we propose an alternative three-step approach to verify whether the global optimum $\mathcal{F}^{*}$ of \eqref{eq:sec3:unknown-proposition} is non-negative. The essence of the min-max problem is to find an input $u \in \mathbb{U}$ that satisfies the DTCBF constraint~\eqref{eq:sec2:DTCBF-constraint-alpha} for each \mbox{$x \in \mathcal{C}$}. In our approach, we determine a suitable control input that satisfies the DTCBF constraint~\eqref{eq:sec2:DTCBF-constraint-alpha} for the state at the center of $\mathbb{X}$ by formulating the inner maximization in \eqref{eq:sec3:unknown-proposition}. If this control input satisfies \eqref{eq:sec2:DTCBF-constraint-alpha} for all states in $\mathbb{X}$, $(h, \gamma)$ is verified as a valid DTCBF. Otherwise, we proceed to partition $\mathbb{X}$ into smaller subdomains and repeat the same procedure within each subdomain. We continue this process until either $(h, \gamma)$ is verified as a valid DTCBF within each subdomain that intersects $\mathcal{C}$, after which a piecewise constant friend $\pi$ is obtained, or a counterexample is found showing that $(h, \gamma)$ is not a DTCBF (within predefined tolerances). 

To simplify the introduction of the algorithm, we initially focus on the verification of a candidate DTCBF within \mbox{a \textit{subdomain}} $X^{(k)} \defeq [x^{lb, (k)}, ~ x^{ub, (k)}] \subseteq \mathbb{X}$, where $k \in \mathbb{N}_0$ is the subdomain number.
\subsubsection{Known Control Policy -- Verification within \texorpdfstring{$X^{(k)}$}{X(k)}} \label{sec:sec4:verification-known-within-subdomain} As discussed in Proposition \ref{prop:sec3:verification-known}, $(h, \gamma)$ with a control policy \mbox{$\pi:\mathcal{C} \rightarrow \mathbb{U}$} is a valid DTCBF within $X^{(k)} \defeq [x^{lb, (k)}, ~ x^{ub, (k)}]$ if \mbox{$\mathcal{F}^{*,(k)} \geqslant 0$}, where $\mathcal{F}^{*,(k)} \in \mathbb{R}$ is the global minimum of
\begin{subequations} \label{eq:sec3:verification:known-general}
    \begin{align} 
    \mathcal{F}^{*,(k)} \defeq \min_{x \in X^{(k)}} ~&~ \underbrace{h(f(x, \pi(x))) - h(x) + \gamma(h(x))}_{\defeq \mathcal{F}(x)} \label{eq:sec3:verification:known-general-objective}\\
    \mathst ~&~ \mathcal{H}(x) \defeq -h(x) \leqslant 0. \label{eq:sec3:verification:known-general-constraint}
    \end{align}
\end{subequations}
Since the optimization problem \eqref{eq:sec3:verification:known-general} is generally non-convex, we aim to compute a tractable lower bound on $\mathcal{F}^{*,(k)}$. To do so, we construct a convex underestimator $\widebreve{\mathcal{F}}^{(k)}:\mathbb{R}^n \rightarrow \mathbb{R}$ of the objective function \eqref{eq:sec3:verification:known-general-objective} and a convex underestimator \mbox{$\widebreve{\mathcal{H}}^{(k)}:\mathbb{R}^n \rightarrow \mathbb{R}$} of the constraint \eqref{eq:sec3:verification:known-general-constraint} on $X^{(k)}$ as
\begin{align} 
    \widebreve{\mathcal{F}}^{(k)}(x) &\defeq \mathcal{F}(x) +\sum_{i=1}^{n}\alpha^{(k)}_{\mathcal{F},i}\bigl(x_i^{lb,(k)} - x_i\bigr)\bigl(x_i^{ub,(k)} - x_i\bigr), \label{eq:sec3:verification:known-objective-convexified} \\
    \widebreve{\mathcal{H}}^{(k)}(x) &\defeq \mathcal{H}(x) +\sum_{i=1}^{n}\alpha^{(k)}_{\mathcal{H},i}\bigl(x_i^{lb,(k)} - x_i\bigr)\bigl(x_i^{ub,(k)} - x_i\bigr), \label{eq:sec3:verification:known-constraint-convexified}
\end{align}
where $\alpha^{(k)}_{\mathcal{F},i}, \alpha^{(k)}_{\mathcal{H},i} \in \mathbb{R}_{\geqslant 0}$, $i \in \{1,2, \hdots, n\}$, are sufficiently large and computed based on one of the methods in \cite{ADJIMAN1998a}. Thus, the convex optimization problem is constructed as
\begin{subequations} \label{eq:sec3:verification:known-convex}
    \begin{align} 
    \widebreve{\mathcal{F}}^{*,(k)} \defeq \min_{x \in X^{(k)}} ~&~ \widebreve{\mathcal{F}}^{(k)}(x) \label{eq:sec3:verification:known-convex-objective}\\
    \mathst ~&~ \widebreve{\mathcal{H}}^{(k)}(x) \leqslant 0. \label{eq:sec3:verification:known-convex-constraint}
    \end{align}
\end{subequations}
Given that for all $x \in X^{(k)}$,
\begin{align}
    \widebreve{\mathcal{F}}^{(k)}(x) &\leqslant \mathcal{F}(x), \label{eq:sec3:verification:known:objective-ineq} \\
   \widebreve{\mathcal{H}}^{(k)}(x) &\leqslant \mathcal{H}(x), \label{eq:sec3:verification:known:constraint-ineq}
\end{align}
it follows that $\mathcal{F}^{*,(k)} \geqslant \widebreve{\mathcal{F}}^{*,(k)}$. 
As a result, by computing the global minimum $\widebreve{\mathcal{F}}^{*,(k)}$ of the convex problem~(11), we obtain a valid lower bound on $\mathcal{F}^{*,(k)}$, which can be used to assess the non-negativity of $\mathcal{F}^{*,(k)}$. In particular, if \eqref{eq:sec3:verification:known-convex} is feasible and $\widebreve{\mathcal{F}}^{*,(k)} \geqslant 0$, then $(h, \gamma)$ with its friend $\pi$ is verified as a valid DTCBF within $X^{(k)}$, since \mbox{$\mathcal{F}^{*,(k)} \geqslant \widebreve{\mathcal{F}}^{*,(k)} \geqslant 0$}. If \eqref{eq:sec3:verification:known-convex} is infeasible, we also approve $X^{(k)}$, since this implies that \eqref{eq:sec3:verification:known-general} is infeasible as well due to \eqref{eq:sec3:verification:known:constraint-ineq}. Otherwise, as explained in ``3) \nameref{sec:sec3:proposed-algorithm:algorithm-overview}'' below, we divide $X^{(k)}$ into smaller subdomains to obtain tighter convex underestimators based on Proposition~\ref{prop:sec3:preliminaries:tighter-underestimator}.
\subsubsection{Unknown Control Policy -- Verification within \texorpdfstring{$X^{(k)}$}{X(k)}}\label{sec:sec4:verification-unknown-within-subdomain} For the unknown control policy case, the three-step approach to verify a candidate DTCBF $h$ with a \mbox{$\gamma \in \mathcal{K}^{\leqslant \mathrm{id}}_{\infty}$} within subdomain $X^{(k)}$ for the \mbox{system \eqref{eq:sec2:dynamical-system}} with $\mathbb{U}$ is proposed as follows:

\paragraph*{Step I} We select the state $\bar{x}^{(k)} \in \mathbb{R}^n$ as the center of the subdomain $X^{(k)} \defeq [x^{lb, (k)}, ~ x^{ub, (k)}]$, i.e., 
\begin{align} \label{eq:sec3:verification:unknown:step1}
    \bar{x}^{(k)}_i \defeq \frac{x^{ub, (k)}_i + x^{lb, (k)}_i}{2}, 
\end{align}
where $i \in \{1, \hdots, n\}$ and $\bar{x}^{(k)} \defeq [\bar{x}^{(k)}_1~ \hdots~ \bar{x}^{(k)}_n ]^{\T}$.

\paragraph*{Step II} \label{algorithm:verification:unknown:step2} For the state $\bar{x}^{(k)}$, we aim to find an admissible control input $u^{*,(k)} \in \mathbb{U}$ that satisfies \mbox{the DTCBF constraint \eqref{eq:sec2:DTCBF-constraint-alpha}}. Therefore, we consider the inner maximization of \eqref{eq:sec3:unknown-proposition} and solve it to global optimality using the $\alpha$BB algorithm \cite{Floudas2010a}: 
\begin{align}
    u^{*,(k)} \in \underset{u \in \mathbb{U}}{\argmax} \hspace{0.03cm} & \hspace{0.03cm} \underbrace{h\bigl(f(\bar{x}^{(k)},u)\bigr) - h\bigl(\bar{x}^{(k)}\bigr) + \gamma\bigl(h(\bar{x}^{(k)})\bigr)}_{\defeq \mathcal{F}^{(k)}_u(u)}. \label{eq:sec3:verification:unknown:step2-objective}
\end{align}

\begin{remark}  
By Assumptions \ref{assump:sec3:C-bounded} and \ref{assump:sec3:continuous-functions}, $\mathbb{U}$ is compact, and $h$ and $f$ are continuous. Thus, based on the Weierstrass theorem \cite{ARORA2017}, a global minimizer $u^{*,(k)} \in \mathbb{U}$ of \eqref{eq:sec3:verification:unknown:step2-objective} exists.
\end{remark}

\paragraph*{Step III} \label{algorithm:verification:unknown:step3} We check whether the control input $u^{*,(k)}$ can satisfy the DTCBF constraint \eqref{eq:sec2:DTCBF-constraint-alpha} for all $x \in X^{(k)}$. This involves verifying whether $\mathcal{F}^{*,(k)} \geqslant 0$, where
\begin{subequations} \label{eq:sec3:verification:unknown:step3-general}
    \begin{align} 
        \mathcal{F}^{*,(k)} \defeq \min_{x \in X^{(k)}} &~ 
        \underbrace{\mystrut{1ex} h(f(x,u^{*,(k)})) - h(x) +\gamma(h(x))}_{\defeq \mathcal{F}(x)}\label{eq:sec3:verification:unknown:step3-general:objective}\\
        \mathst ~&~ \mathcal{H}(x) \defeq -h(x) \leqslant 0. \label{eq:sec3:verification:unknown:step3-general:constraint}
    \end{align}
\end{subequations}
Similarly, since \eqref{eq:sec3:verification:unknown:step3-general} is generally non-convex, we construct the convex optimization problem as%
\begin{subequations}\label{eq:sec3:verification:unknown:step3-convex}
    \begin{align} 
        \widebreve{\mathcal{F}}^{*,(k)}  \defeq \min_{x \in X^{(k)}}  ~&~ \widebreve{\mathcal{F}}^{(k)}(x) \label{eq:sec3:verification:unknown:step3-convex-objective} \\
        \mathst ~&~ \widebreve{\mathcal{H}}^{(k)}(x) \leqslant 0, \label{eq:sec3:verification:unknown:step3-convex-constraint}
    \end{align}
\end{subequations}
where $\widebreve{\mathcal{F}}^{(k)}$ and $\widebreve{\mathcal{H}}^{(k)}$ are convex underestimators of the objective function \eqref{eq:sec3:verification:unknown:step3-general:objective} and the constraint \eqref{eq:sec3:verification:unknown:step3-general:constraint} on $X^{(k)}$, constructed according to \eqref{eq:sec3:verification:known-objective-convexified} and \eqref{eq:sec3:verification:known-constraint-convexified}, respectively. Similar to the known control policy case, $(h, \gamma)$, with its friend \mbox{$\pi(x) \defeq u^{*,(k)}$} for $x\in X^{(k)}$, is verified as a valid DTCBF within $X^{(k)}$ if $\widebreve{\mathcal{F}}^{*,(k)} \geqslant 0$.
\subsubsection{Algorithm Overview} \label{sec:sec3:proposed-algorithm:algorithm-overview}
By solving the optimization problem \eqref{eq:sec3:verification:known-convex} for the known control policy case or \eqref{eq:sec3:verification:unknown:step3-convex} for the unknown control policy case on a subdomain $X^{(k)}$, we encounter three cases:
\begin{enumerate}[label=(\Alph*), ref=(\Alph*)]
    \item  $\widebreve{\mathcal{F}}^{*,(k)} \geqslant 0$: $(h, \gamma)$ is verified as a valid DTCBF within $X^{(k)}$ as \mbox{$\mathcal{F}^{*,(k)} \geqslant \widebreve{\mathcal{F}}^{*,(k)} \geqslant 0$}. Thus, we approve this subdomain. \label{case:convex-verification-known-A}
    \item The convex optimization problem \eqref{eq:sec3:verification:known-convex} or \eqref{eq:sec3:verification:unknown:step3-convex} is infeasible: Given that \eqref{eq:sec3:verification:known:constraint-ineq} holds, it follows that $X^{(k)} \cap \mathcal{C} = \emptyset$. Thus, we approve this subdomain. \label{case:convex-verification-known-B}
    \item \label{case:convex-verification-known-C} $\widebreve{\mathcal{F}}^{*,(k)} < 0$: For the known control policy case, consider $x^{*,(k)}$ as the global minimizer of \eqref{eq:sec3:verification:known-convex} and $\mathcal{F}$ as defined in \eqref{eq:sec3:verification:known-general}.
    For the unknown control policy case, consider $\bar{x}^{(k)}$ as obtained from \eqref{eq:sec3:verification:unknown:step1}, $\mathcal{F}^{(k)}_u$ as defined in \eqref{eq:sec3:verification:unknown:step2-objective}, and \mbox{$u^{*,(k)}$} as a global minimizer of \eqref{eq:sec3:verification:unknown:step2-objective}.
    \begin{enumerate} [label=(C.\arabic*), ref=(C.\arabic*), leftmargin=18pt]
        \item $x^{*,(k)} \in \mathcal{C}$ and $\mathcal{F}(x^{*,(k)}) < 0$ for the known control policy, or $\bar{x}^{(k)} \in \mathcal{C}$, $u^{*,(k)} \in \mathbb{U}$, and \mbox{$\mathcal{F}^{(k)}_u(u^{*,(k)}) < 0$} for the unknown control policy: We terminate the algorithm and report that $(h, \gamma)$ is not a valid DTCBF, or that $\pi$, if given, is not a friend \textit{(without any conservatism)}, providing $x^{*,(k)}$ or $\bar{x}^{(k)}$ as a counterexample. \label{case:convex-verification-known-C-1}
        \item Otherwise: $X^{(k)}$ is divided into smaller subdomains. \label{case:convex-verification-known-C-2}
    \end{enumerate}
\end{enumerate}

Moreover, to ensure that the algorithm terminates in a finite number of iterations, certain stopping criteria are imposed. In the case of a known control policy, we impose stopping criteria on the maximum separation between the objective function \eqref{eq:sec3:verification:known-general-objective} and its convex underestimator \eqref{eq:sec3:verification:known-convex-objective}, as well as the constraint \eqref{eq:sec3:verification:known-general-constraint} and its convex underestimator \eqref{eq:sec3:verification:known-convex-constraint}. Specifically, the algorithm terminates on a subdomain $X^{(k)}$ if
\begin{align}
    \frac{\max_i\alpha_{\mathcal{F},i}^{(k)}}{4}\sum_{i=1}^{n} \bigl(x_i^{ub,(k)} - x_i^{lb,(k)}\bigr)^2 &\leqslant \epsilon_f \label{eq:sec3:verification:stopping-criteria-1},\\
    \frac{\max_i\alpha_{\mathcal{H},i}^{(k)}}{4}\sum_{i=1}^{n}\bigl(x_i^{ub,(k)} - x_i^{lb,(k)}\bigr)^2 &\leqslant \epsilon_h \label{eq:sec3:verification:stopping-criteria-2}, 
\end{align}
where $\epsilon_f, \epsilon_h \in \mathbb{R}_{>0}$ are predefined values. In the case of an unknown control policy, we impose the same stopping criteria but additionally consider the size of subdomains, which determines the smallest domain on which the resulting control policy $\pi$ can assign a distinct control input. In particular, the algorithm terminates on $X^{(k)}$ if \eqref{eq:sec3:verification:stopping-criteria-1}, \eqref{eq:sec3:verification:stopping-criteria-2}, and 
\begin{align}
    \delta_k^2 \defeq \sum_{i=1}^{n}\bigl(x_i^{ub,(k)} - x_i^{lb,(k)}\bigr)^2 &\leqslant \epsilon_d \label{eq:sec3:verification:stopping-criteria-3}
\end{align}
are met, where $\delta_k \in \mathbb{R}_{> 0}$ is the diagonal size of $X^{(k)}$ and \mbox{$\epsilon_d \in \mathbb{R}_{>0}$} is a predefined value. In this case, a piecewise constant control policy $\pi$ is computed, which can assign a distinct control input $u \in \mathbb{U}$ to each subdomain $X^{(k)} \subseteq \mathbb{X}$ whose diagonal size $\delta_k$ is larger than $\sqrt{\epsilon_d}$. If the stopping criteria are met, we return the state $\bar{x}^{(k)}$ at the center of $X^{(k)}$ for both cases as a counterexample to \textit{the selected conservatism}.
The overall BB algorithm is detailed in Algorithm~\ref{algorithm:verification}. 
\begin{remark}
If the stopping criteria are satisfied, the algorithm may fail to verify a valid DTCBF (though it never verifies a non-valid one). As proved in Theorem \ref{theorem:verification-conservatism} below, choosing lower values for $\epsilon_f, \epsilon_h, \epsilon_d \in \mathbb{R}_{> 0}$ in the stopping criteria increases the accuracy of the verification algorithm (in terms of not failing to verify a valid DTCBF), but leads to longer computation times due to the need for finer refinements based on Theorem \ref{theorem:verification-finite} below. Thus, the choice of these values depends on the application and involves a trade-off between computation time and accuracy of the outcome.
\end{remark}
\begin{algorithm}[t]
\caption{Algorithm for Verification} \label{algorithm:verification}
\begin{algorithmic}[1] 
\State \textbf{Input:} $f$, $h$, $\gamma$, $\pi$ (if given). \mbox{$\epsilon_f, \epsilon_h \in \mathbb{R}_{>0}$}, and additionally $\epsilon_d \in \mathbb{R}_{>0}$ if the control policy is unknown. An $n$-rectangle set $\mathbb{X}$ such that $\mathcal{C}\subseteq \mathbb{X}$.
\State $X^{(0)} \gets \mathbb{X}$
\State $\mathcal{L} \gets \{0\}$ \Comment{$\mathcal{L}$ is the list of subdomain numbers that remain to be handled}
\State $n_{dom} \gets 0$ \Comment{Number of subdomains}
\While{$\mathcal{L} \neq \emptyset$} \label{algorithm:while}
\State $k \gets \text{any subdomain number from the list } \mathcal{L}$
\If{control policy is known}
    \State solve the optimization problem \eqref{eq:sec3:verification:known-convex} on $X^{(k)}$
\ElsIf{control policy is unknown}
        \State do the three steps in Section \ref{sec:sec4:verification-unknown-within-subdomain} for $X^{(k)}$
\EndIf
\If{Case \ref{case:convex-verification-known-A} or Case \ref{case:convex-verification-known-B}}
    \State $\mathcal{L} \gets \mathcal{L}\backslash \{k\}$
\ElsIf{Case \ref{case:convex-verification-known-C-1}}
    \State \textbf{print} ``$(h,\gamma)$ is not a valid DTCBF {\small$[$or $\pi$, if given, is not a friend$]$}, and the counterexample is:'' $\bar{x}^{(k)}$ {\small$[$or $x^{*,(k)}$$]$}
    \State terminate the algorithm
\Else \Comment{Case \ref{case:convex-verification-known-C-2}}
\If{stopping criteria are met on $X^{(k)}$}
    \State \textbf{print} ``$(h,\gamma)$ is not valid {\small $[$or $\pi$, if given, is not a friend$]$} \textit{for the selected conservatism} $\epsilon_f, \epsilon_h$ {\small$[$and $\epsilon_d]$}. The counterexample is:'' $\bar{x}^{(k)}$
    \State terminate the algorithm
\EndIf
\State divide $X^{(k)}$ into $X^{(n_{dom}+1)}$ and $X^{(n_{dom}+2)}$ using the \textit{scaled longest side strategy} \cite{Maranas1994b}
\State $\mathcal{L} \gets (\mathcal{L}\backslash\{k\}) \cup \{n_{dom}+1,~ n_{dom}+2\}$
\State $n_{dom} \gets n_{dom}+2$
\EndIf
\EndWhile
\State \textbf{print} ``$(h,\gamma)$ is a DTCBF"
\end{algorithmic}
\end{algorithm}
\begin{remark} \label{remark:differences-ourmethod-abb}
    The main differences between our method and the $\alpha$BB algorithm \cite{Floudas2010a} are as follows:
    \begin{itemize}
        \item For the known control policy case, we compute only lower bounds on the global minimum of \eqref{eq:sec3:known-proposition} (through the optimization problem \eqref{eq:sec3:verification:known-convex}), not upper bounds. We also approve subdomains $X^{(k)}$ where $\widebreve{\mathcal{F}}^{*,(k)} \geqslant 0$. Thus, our method can be significantly faster than the $\alpha$BB algorithm for the verification problem. Moreover, it can falsify invalid candidate DTCBFs by providing counterexamples without conservatism, which is not possible with the $\alpha$BB algorithm due to the concept of $\epsilon$-feasibility (see Section~\ref{sec:sec5:verification}). 
        \item For the unknown control policy case, the $\alpha$BB algorithm is incapable of directly solving \eqref{eq:sec3:unknown-proposition} and thus cannot verify or falsify a candidate DTCBF.
    \end{itemize}
\end{remark}
\begin{theorem} \label{theorem:verification-finite}
    Given $\epsilon_f, \epsilon_h \in \mathbb{R}_{> 0}$ for the known control policy case or $\epsilon_f, \epsilon_h, \epsilon_d \in \mathbb{R}_{> 0}$ for the unknown control policy case, Algorithm \ref{algorithm:verification} terminates in a finite number of iterations.
    
    \begin{proof}
      The algorithm terminates either if $\mathcal{L} = \emptyset$ or the stopping criteria, defined by \eqref{eq:sec3:verification:stopping-criteria-1} and \eqref{eq:sec3:verification:stopping-criteria-2} for the known control policy case or \eqref{eq:sec3:verification:stopping-criteria-1}\textendash \eqref{eq:sec3:verification:stopping-criteria-3} for the unknown control policy case, are met on a subdomain. Assuming $\max_i\alpha^{(k)}_{\mathcal{F}, i} \neq 0$ on a subdomain $X^{(k)}$, we rewrite the stopping criterion \eqref{eq:sec3:verification:stopping-criteria-1} as
        \begin{align} 
        \delta_k^2 \defeq \sum_{i=1}^{n} \left(x_i^{ub,(k)} - x_i^{lb,(k)}\right)^2 &\leqslant \frac{4\epsilon_f}{\max_i\alpha^{(k)}_{\mathcal{F},i}}. \nonumber
        \end{align}
        With the same reasoning for \eqref{eq:sec3:verification:stopping-criteria-2} assuming $\max_i\alpha^{(k)}_{\mathcal{H},i} \neq 0$, we conclude that the algorithm terminates for the known control policy case if the diagonal size $\delta_k$ of $X^{(k)}$ becomes
        \begin{align} \label{eq:sec3:verification:proof-finite-ineq}
            \delta_k \leqslant \mathmin \biggl\{ \sqrt{\frac{4\epsilon_f}{\max_i\alpha^{(k)}_{\mathcal{F},i}}}, ~ \sqrt{\frac{4\epsilon_h}{\max_i\alpha^{(k)}_{\mathcal{H},i}}} \biggr\}.
        \end{align}
        Similarly, in the case where a corresponding control policy is unknown, the algorithm terminates if the diagonal size of $X^{(k)}$ becomes
        \begin{align} \label{eq:sec3:verification:proof-finite-ineq-2}
            \delta_k \leqslant \mathmin \biggl\{ \sqrt{\frac{4\epsilon_f}{\max_i\alpha^{(k)}_{\mathcal{F},i}}}, ~ \sqrt{\frac{4\epsilon_h}{\max_i\alpha^{(k)}_{\mathcal{H},i}}}, ~ \sqrt{\epsilon_d} \biggr\}.
        \end{align}
        Based on Proposition~\ref{prop:sec3:preliminaries:tighter-underestimator}, it holds that $\alpha^{(k+1)}_{\mathcal{F},i} \leqslant \alpha^{(k)}_{\mathcal{F},i}$ and $ \alpha^{(k+1)}_{\mathcal{H},i} \leqslant \alpha^{(k)}_{\mathcal{H},i}$ for all $k \in \mathbb{N}_0$, as $X^{(k+1)} \subset X^{(k)}$. Moreover, $\epsilon_f, \epsilon_h, \epsilon_d \in \mathbb{R}_{> 0}$ are fixed constants. Thus, the right-hand sides of \eqref{eq:sec3:verification:proof-finite-ineq} and \eqref{eq:sec3:verification:proof-finite-ineq-2} are non-zero and non-decreasing as $k$ increases, and as \cite{Maranas1994b} proves, the number of dividing iterations required using the scaled longest side strategy to obtain a subdomain $X^{(k)}$ with diagonal size of $\delta_{k} \in \mathbb{R}_{>0}$ is finite, so Algorithm \ref{algorithm:verification} terminates in a finite number of iterations. If $\max_i\alpha^{(k)}_{\mathcal{F}, i} = 0$ or $\max_i\alpha^{(k)}_{\mathcal{H}, i} = 0$, the corresponding stopping criterion is met, regardless of the subdomain size. 
    \end{proof}
\end{theorem}
\begin{remark}
    Note that $\max_i\alpha^{(k)}_{\mathcal{F}, i} = 0$ and $\max_i\alpha^{(k)}_{\mathcal{H}, i} = 0$ imply that $\mathcal{F}$ and $\mathcal{H}$ are convex, respectively, on $X^{(k)}$. In this case, for a known control policy, since \eqref{eq:sec3:verification:known-convex} is equivalent to \eqref{eq:sec3:verification:known-general}, if $\widebreve{\mathcal{F}}^{*,(k)} < 0$, the global minimizer of \eqref{eq:sec3:verification:known-convex} already serves as a counterexample, meaning Case~\ref{case:convex-verification-known-C-1} is met before the stopping criteria are checked. For an unknown control policy, although \eqref{eq:sec3:verification:unknown:step3-general} and \eqref{eq:sec3:verification:unknown:step3-convex} are equivalent in this case, the algorithm proceeds while $\widebreve{\mathcal{F}}^{*,(k)} < 0$ until the diagonal size of $X^{(k)}$ is at most $\sqrt{\epsilon_d}$, allowing distinct control inputs to be applied. If \mbox{$\widebreve{\mathcal{F}}^{*,(k)} \geqslant 0$}, $X^{(k)}$ is approved for both cases. 
\end{remark}
\begin{theorem} \label{theorem:verification-conservatism}
    Consider the system \eqref{eq:sec2:dynamical-system} with \mbox{$\mathbb{U}$}, a \mbox{$\gamma \in \mathcal{K}^{\leqslant \mathrm{id}}_{\infty}$}, and a candidate DTCBF $h$. Additionally, let \mbox{$\epsilon_f, \epsilon_h \in \mathbb{R}_{>0}$} and a control policy $\pi:\mathcal{C} \rightarrow \mathbb{U}$ be given for the known control policy case, or let $\epsilon_f, \epsilon_h, \epsilon_d \in \mathbb{R}_{>0}$ be given for the unknown control policy case. In the case of a known control policy $\pi$, Algorithm~\ref{algorithm:verification} either verifies that $h$ with $\gamma$ satisfies 
    \begin{align} \label{eq:sec3:verification:proof-DTCBF}
        h(f(x,\pi(x))) - h(x) + \gamma(h(x)) \geqslant \epsilon_f, 
    \end{align}
    for all $x \in \{\xi \in \mathbb{R}^n \mid h(\xi)  + \epsilon_h \geqslant 0\}$, or falsifies it by providing a counterexample. For the unknown control policy case, Algorithm \ref{algorithm:verification} either provides a control policy $\pi:\mathcal{C} \rightarrow \mathbb{U}$ such that \eqref{eq:sec3:verification:proof-DTCBF} is satisfied for all $x \in \{\xi \in \mathbb{R}^n \mid h(\xi)  + \epsilon_h \geqslant 0\}$, where $\pi$ is a piecewise constant friend that can apply a distinct $u \in \mathbb{U}$ for each subdomain with diagonal size larger than $\sqrt{\epsilon_d}$, or falsifies it by providing a counterexample.
    
    \begin{proof}
        Before giving the proof, it should be noted that verifying whether \eqref{eq:sec3:verification:proof-DTCBF} is satisfied implies determining whether $(h,\gamma)$ with its friend $\pi$ is a DTCBF for the given conservatism.
        
        Assume that the stopping criteria are met for a subdomain $X^{(k)}$. Based on Proposition \ref{prop:sec3:preliminaries:max-separation} and the stopping criteria, the maximum separation between the objective function \eqref{eq:sec3:verification:known-general-objective} or \eqref{eq:sec3:verification:unknown:step3-general:objective} and its convex underestimator \eqref{eq:sec3:verification:known-convex-objective} or \eqref{eq:sec3:verification:unknown:step3-convex-objective} on $X^{(k)}$ occurs at the center of $X^{(k)}$ and is upper bounded by
        \begin{align}
            &\frac{1}{4}\sum_{i=1}^{n} \alpha^{(k)}_{\mathcal{F}, i}\bigl(x_i^{ub,(k)} - x_i^{lb,(k)}\bigr)^2  \leqslant \nonumber \\
            & \qquad \qquad  \frac{\max_i\alpha^{(k)}_{\mathcal{F}, i}}{4}\sum_{i=1}^{n} \bigl(x_i^{ub,(k)} - x_i^{lb,(k)}\bigr)^2 \leqslant  \epsilon_f. \nonumber
        \end{align}
        Similarly, the maximum separation between the constraint \eqref{eq:sec3:verification:known-general-constraint} or \eqref{eq:sec3:verification:unknown:step3-general:constraint} and its convex underestimator \eqref{eq:sec3:verification:known-convex-constraint} or \eqref{eq:sec3:verification:unknown:step3-convex-constraint} on $X^{(k)}$ is upper bounded by 
        \begin{align}
            \frac{1}{4}\sum_{i=1}^{n} \alpha^{(k)}_{\mathcal{H}, i}\bigl(x_i^{ub,(k)} - x_i^{lb,(k)}\bigr)^2 &\leqslant \epsilon_h. \nonumber
        \end{align}
        Based on Algorithm \ref{algorithm:verification}, $k \in \mathcal{L}$ implies that $\widebreve{\mathcal{F}}^{*,(k)} < 0$. Thus,
        \begin{align}
            h(f(\bar{x}^{(k)},\bar{u})) - h(\bar{x}^{(k)}) + \gamma(h(\bar{x}^{(k)})) < \epsilon_f, \nonumber
        \end{align}
        where $\bar{x}^{(k)} \in \mathbb{R}^n$ is the state at the center of $X^{(k)}$, and \mbox{$h(\bar{x}^{(k)}) \geqslant - \epsilon_h$}. If the control policy is known, $\bar{u} \defeq \pi(\bar{x}^{(k)})$. If the control policy is not \textit{a priori} known, as indicated in Step~II, an input \mbox{$\bar{u} \defeq u^{*,(k)} \in \mathbb{U}$} is computed for $\bar{x}^{(k)}$.
        
        On the other hand, if $\mathcal{L} = \emptyset$, it implies that \eqref{eq:sec3:verification:proof-DTCBF} is satisfied. Thus, $(h,\gamma)$ with its friend $\pi$ is a DTCBF.
    \end{proof}
\end{theorem}

\section{Synthesis}
\label{sec:sec4}
In this section, we consider the problem of synthesizing DTCBFs.

\subsection{Problem Statement and Assumptions} \label{sec:sec4:problem-statement} 
\begin{problem}[Synthesis] \label{problem:synthesis}
    Consider the system \eqref{eq:sec2:dynamical-system} with the control admissible set $\mathbb{U}$, the safe set $\mathcal{S}$ in \eqref{eq:sec2:safe-set}, and a given \mbox{$\gamma \in \mathcal{K}^{\leqslant \mathrm{id}}_{\infty}$}. The objective is to synthesize a function \mbox{$h:\mathbb{R}^n \rightarrow \mathbb{R}$} with zero-superlevel set $\mathcal{C}$ and a control policy~$\pi$ such that the size of $\mathcal{C}$ is maximized (in an appropriate sense), and the following properties are satisfied:
    \begin{enumerate}[ref =\arabic{enumi}, label=  \arabic{enumi})]
        \item $(h, \gamma)$ with the control policy $\pi$ satisfies the DTCBF constraint \eqref{eq:sec2:DTCBF-constraint-alpha} for all $x \in \mathcal{C}$. \label{condition:synthesis-1}
        \item $\pi$ is admissible for all $x \in \mathcal{C}$ with respect to $\mathbb{U}$, and thus, together with Property \ref{condition:synthesis-1}, $\pi$ is a friend of $(h, \gamma)$. \label{condition:synthesis-2}
        \item $\mathcal{C}$ is a subset of the safe set, i.e., $\mathcal{C} \subseteq \mathcal{S}$. Thus, the system \eqref{eq:sec2:dynamical-system} with $\mathbb{U}$ is $(\mathcal{S}, X_0)$-safe for some $X_0 \subseteq \mathcal{C}$. \label{condition:synthesis-3}
    \end{enumerate}
    \end{problem}
    
    To address the synthesis problem in a tractable manner, we parameterize $h$ as
    \begin{align} \label{eq:sec4:parameterized-DTCBF}
        h(x;\vartheta) \defeq \sum_{i = 1}^{n_\vartheta} \Lambda_i(x;\vartheta_i),
    \end{align}
    where $\Lambda_i:\mathbb{R}^{n} \times \mathbb{R}\rightarrow \mathbb{R}$, $i \in \{1,2, \hdots, n_{\vartheta}\}$, is a given basis function, $\vartheta \defeq [\vartheta_1 ~\hdots ~\vartheta_{n_{\vartheta}}]^{\T} \in \Theta \subseteq \mathbb{R}^{n_{\vartheta}}$ is the vector of unknown coefficients, and \mbox{$n_\vartheta \in \mathbb{N}$} is the number of coefficients. With $h$ parameterized by $\vartheta$, we refer to its zero-superlevel set as $\mathcal{C}(\vartheta)$. Likewise, we parameterize the control policy $\pi$ as
    \begin{align} \label{eq:sec4:parameterized-control-policy}
        \pi_i(x;\mu_i) \defeq \sum_{j = 1}^{n_{\mu_i}} \Pi_{i,j}(x;\mu_{i,j}), ~~ i \in \{1,2,\hdots, m\}, 
    \end{align}
    where $\pi(x; \mu) \defeq [\pi_1(x;\mu_1)~\hdots~ \pi_m(x; \mu_m)]^{\T}$. Here, \mbox{$\Pi_{i,j}: \mathbb{R}^{n} \times \mathbb{R} \rightarrow \mathbb{R}$} is a given basis function, $\mu_i \defeq [\mu_{i,1} ~\hdots ~\mu_{i,n_{\mu_i}}]^{\T} \in \mathbb{R}^{n_{\mu_i}}$ is the vector of unknown coefficients, and $n_{\mu_i} \in \mathbb{N}$ is the number of coefficients in $\pi_i$, \mbox{$i \in \{1,2,\hdots, m\}$}. Moreover, \mbox{$\mu \defeq [\mu_1~ \hdots~\mu_m]^{\T} \in M \subseteq \mathbb{R}^{n_\mu}$} is the vector of all unknown coefficients in $\pi$, and \mbox{$n_\mu \defeq \sum_{i=1}^{m} n_{\mu_i}$}
    is the total number of unknown coefficients. 
    We emphasize that the proposed bilevel approach is applicable to any parameterization of $h$ and $\pi$, including those that are nonlinear in the coefficients.
\begin{remark}
    Regarding the choice of basis functions in $h$ and $\pi$, using richer basis functions increases the likelihood of successfully synthesizing a valid DTCBF but also leads to a higher (offline) computational burden for synthesis. Therefore, a trade-off must be considered. However, it is advisable to avoid complicated basis functions for $h$, as the synthesized DTCBF is used in an online optimization problem of the form%
    \begin{subequations} \label{eq:remark-online}
    \begin{align}
        u_t^* \in \text{argmin}_{u \in \mathbb{U}} ~&~ \lVert u -  \pi_{\text{nom}}(x_t)\rVert \\
        \text{s.t.}~&~ h(f(x_t,u)) -h(x_t) \geqslant - \gamma(h(x_t))
    \end{align}
    \end{subequations}
    to compute a safe control input $u_t^* \in \mathbb{U}$ at each time instant $t \in \mathbb{N}_0$ for the current state $x_t$, given a nominal controller \mbox{$\pi_{\text{nom}}: \mathbb{R}^n \rightarrow \mathbb{R}^m$}. Thus, quadratic basis functions for both $h$ and $\pi$ are often used in practice. If a DTCBF cannot be obtained using quadratic basis functions for a specific system, the basis functions can be enriched to improve expressiveness, at the cost of higher offline and online computational complexity.
\end{remark}
\begin{remark} \label{sec4:remark-existence-class-kappa}
    In Problem \ref{problem:synthesis}, we assume, without loss of generality, that $\gamma \in \mathcal{K}^{\leqslant \mathrm{id}}_{\infty}$ is given. If $\gamma$ is not specified, we can parameterize it similarly, e.g., by choosing $\gamma(x) \defeq cx$, where $c \in (0,1] \subset \mathbb{R}$ is unknown. The proposed bilevel optimization approach can then be used to synthesize both $h$ and $\gamma$.
\end{remark}
\begin{assumption} \label{assump:sec4:bounded-DTCBF}
    $\mathcal{C}(\vartheta)$ is compact for all \mbox{$\vartheta \in \Theta$}.
\end{assumption}
\begin{assumption} \label{assump:sec4:compact-unknown-coefficients}
    $\Theta \subset \mathbb{R}^{n_\vartheta}$ and $M \subset \mathbb{R}^{n_\mu}$ are both compact. 
\end{assumption}
\begin{assumption} \label{assump:sec4:twice-differentiable}
    The mapping $f$ associated with the system~\eqref{eq:sec2:dynamical-system}, the mapping $s$ defining the safe set $\mathcal{S}$ in \eqref{eq:sec2:safe-set}, \mbox{$\gamma \in \mathcal{K}^{\leqslant \mathrm{id}}_{\infty}$}, the basis functions $\Lambda_i$, $i \in \{1,2,\hdots, n_{\vartheta}\}$, and the basis functions $\Pi_{i,j}$, \mbox{$i \in \{1,2, \hdots, m\}$} and \mbox{$j \in \{1,2, \hdots, n_{\mu_i}\}$}, are all twice continuously differentiable.
\end{assumption}
\begin{assumption} \label{assump:sec4:chosen-X}
    An $n$-rectangle set $\mathbb{X}$ is chosen such that $\mathcal{C}(\vartheta) \subset \mathbb{X}$ and $\partial \mathcal{C}(\vartheta) \cap \partial \mathbb{X} = \emptyset$  hold for all $\vartheta \in \Theta$.
\end{assumption}

The proposed bilevel approach is applicable to common structures of the control admissible set $\mathbb{U} \subseteq \mathbb{R}^m$, such as semi-algebraic sets and $m$-rectangle sets. In this paper, we consider Assumption \ref{assump:sec4:m-rec-control-admissible} as a specific case to explain our method.
\begin{assumption} \label{assump:sec4:m-rec-control-admissible}
    $\mathbb{U}$ is an \mbox{$m$-rectangle set}.
\end{assumption}
\begin{assumption} \label{assump:sec4:MFCQ}
    The function $h$ is parameterized such that $\nabla_x h(x; \vartheta) \neq [0~ \hdots~ 0]^{\T} \in \mathbb{R}^n$ for all $\vartheta \in \Theta$ and for all $x \in \mathbb{X}$ such that $h(x; \vartheta) = 0$.  
\end{assumption}
Assumption~\ref{assump:sec4:MFCQ} holds for commonly used DTCBFs, including quadratic ones. 
\subsection{Bilevel Approach for Synthesis} \label{sec:sec4:optimization-framework}
We continue by formalizing the synthesis problem in a numerical optimization framework. To this end, each of the three properties in Problem \ref{problem:synthesis} is formulated as a subproblem of the overarching synthesis problem.

Property \ref{condition:synthesis-1} is closely related to the verification problem (Problem \ref{problem:verification}) as discussed in Section \ref{sec:sec3}. To recall, the function $h$, parameterized by $\vartheta$, with a given $\gamma \in \mathcal{K}^{\leqslant \mathrm{id}}_{\infty}$ and the control policy $\pi$, parameterized by $\mu$, is a DTCBF if and only if
\begin{align} \label{eq:sec4:proposed-alg:verification-DTCBF-ineq}
    \mathcal{H}(\xi_h^*; \vartheta, \mu) &\defeq h\left(f(\xi_h^*,\pi(\xi_h^*;\mu));\vartheta\right) - h(\xi_h^*;\vartheta)  \nonumber \\
    & \quad + \gamma(h(\xi_h^*;\vartheta)) \geqslant 0,
\end{align}
where $\xi_h^* \in \mathbb{R}^n$ is a global minimizer of
\begin{subequations} \label{eq:sec4:proposed-alg:verification-DTCBF}
    \begin{align}
    \xi_h^* \in \underset{ \xi_h \in \mathbb{X} }{\argmin} ~&~ \mathcal{H}(\xi_h;\vartheta, \mu) \\
    \mathst ~&~ -h(\xi_h;\vartheta) \leqslant 0. \label{eq:sec4:proposed-alg:verification-DTCBF-constraint}
    \end{align}
\end{subequations}
Regarding Property \ref{condition:synthesis-2}, the control policy $\pi$, parameterized by $\mu$, is admissible for all $x \in \mathcal{C}(\vartheta)$ with respect to the control admissible set \mbox{$\mathbb{U} \defeq \left[u^{lb}, ~ u^{ub}\right]$}, and consequently, $\pi$ is a friend of $(h, \gamma)$, if and only if $\pi(x; \mu) \in \mathbb{U}$ for all $x \in \mathcal{C}(\vartheta)$ or equivalently 
\begin{align} 
    u_i^{lb} &\leqslant \pi_i(\xi^*_{\underline{\pi}_i};\mu_i), \label{eq:sec4:proposed-alg:admissible-control-input-1-ineq} \\
    \pi_i(\xi^*_{\overline{\pi}_i};\mu_i) &\leqslant u_i^{ub}, \label{eq:sec4:proposed-alg:admissible-control-input-2-ineq}
\end{align}
for all $i \in \{1,2, \hdots, m\}$, where $\xi^*_{\underline{\pi}_i} \in \mathbb{R}^n$ is a global minimizer of
\begin{subequations} \label{eq:sec4:proposed-alg:admissible-control-input-1}
    \begin{align}
     \xi^*_{\underline{\pi}_i} \in \underset{ \xi_{\underline{\pi}_i} \in \mathbb{X} }{\argmin} ~&~ \pi_i(\xi_{\underline{\pi}_i};\mu_i) \\
    \mathst ~&~ -h(\xi_{\underline{\pi}_i}; \vartheta) \leqslant 0, \label{eq:sec4:proposed-alg:admissible-control-input-1-constraint}
    \end{align}
    \end{subequations}
and $\xi^*_{\overline{\pi}_i} \in \mathbb{R}^n$ is a global minimizer of
\begin{subequations} \label{eq:sec4:proposed-alg:admissible-control-input-2}
    \begin{align}
    \xi^*_{\overline{\pi}_i} \in \underset{ \xi_{\overline{\pi}_i} \in \mathbb{X} }{\argmin} ~&~-\pi_i(\xi_{\overline{\pi}_i};\mu_i) \\
    \mathst ~&~ -h(\xi_{\overline{\pi}_i};\vartheta) \leqslant 0. \label{eq:sec4:proposed-alg:admissible-control-input-2-constraint}
    \end{align}
\end{subequations}
Proceeding with Property \ref{condition:synthesis-3}, we require $\mathcal{C}(\vartheta)$ to be a subset of the safe set, which is equivalent to satisfying
\begin{align} \label{eq:sec4:proposed-alg:subset-safe-set-ineq}
    s(\xi^*_s) \geqslant 0,
\end{align}
where $\xi^*_s \in \mathbb{R}^n$ is a global minimizer of
\begin{subequations} \label{eq:sec4:proposed-alg:subset-safe-set}
    \begin{align}
        \xi^*_s \in \underset{ \xi_s \in \mathbb{X} }{\argmin} ~&~ s(\xi_s) \\
        \mathst ~&~ -h(\xi_s;\vartheta) \leqslant 0. \label{eq:sec4:proposed-alg:subset-safe-set-constraint}
    \end{align}
\end{subequations}

As discussed, the optimization problems \eqref{eq:sec4:proposed-alg:verification-DTCBF}, \eqref{eq:sec4:proposed-alg:admissible-control-input-1}, \eqref{eq:sec4:proposed-alg:admissible-control-input-2}, and \eqref{eq:sec4:proposed-alg:subset-safe-set} are employed to verify Properties \ref{condition:synthesis-1}\textendash\ref{condition:synthesis-3} for given values of $\vartheta$ and $\mu$, thus they serve as inner problems within our proposed bilevel optimization framework. In the outer problem, our objective is to determine optimal values of $\vartheta$ and $\mu$ that not only satisfy constraints \eqref{eq:sec4:proposed-alg:verification-DTCBF-ineq}, \eqref{eq:sec4:proposed-alg:admissible-control-input-1-ineq}, \eqref{eq:sec4:proposed-alg:admissible-control-input-2-ineq}, and \eqref{eq:sec4:proposed-alg:subset-safe-set-ineq} but also maximize the size of $\mathcal{C}(\vartheta)$.

As the optimization problems \eqref{eq:sec4:proposed-alg:verification-DTCBF}, \eqref{eq:sec4:proposed-alg:admissible-control-input-1}, \eqref{eq:sec4:proposed-alg:admissible-control-input-2}, and \eqref{eq:sec4:proposed-alg:subset-safe-set} are independent of each other for fixed values of $\vartheta$ and $\mu$, we can combine them into one inner problem with the constraints \eqref{eq:sec4:proposed-alg:verification-DTCBF-constraint}, \eqref{eq:sec4:proposed-alg:admissible-control-input-1-constraint}, \eqref{eq:sec4:proposed-alg:admissible-control-input-2-constraint}, and \eqref{eq:sec4:proposed-alg:subset-safe-set-constraint} and with the objective function 
\begin{align} \label{eq:sec5:proposed-synthesis:objective-function}
    \mathcal{F}(\eta) &\defeq  \mathcal{H}(\xi_{h};\vartheta, \mu)  \nonumber \\
    &\quad +  \sum_{i = 1}^m \left[\pi_i(\xi_{\underline{\pi}_i}; \mu_i) -  \pi_i(\xi_{\overline{\pi}_i}; \mu_i)\right] +  s(\xi_{s}), \nonumber
\end{align}
where $\eta \defeq [\xi_{h} ~\xi_{\underline{\pi}} ~\xi_{\overline{\pi}} ~ \xi_{s} ~ \vartheta~ \mu]^{\T}$, $\xi_{\underline{\pi}} \defeq [\xi_{\underline{\pi}_1}~ \hdots~ \xi_{\underline{\pi}_m}]^{\T}$, and $\xi_{\overline{\pi}} \defeq [\xi_{\overline{\pi}_1}~ \hdots~ \xi_{\overline{\pi}_m}]^{\T}$. As a result, the bilevel optimization problem to address the synthesis problem, as stated in Problem~\ref{problem:synthesis}, is constructed as 
\begin{subequations} \label{eq:sec4:bilevel-synthesis}
	\begin{align}
	\rho^* \hspace{-0.03cm} \defeq	\hspace{-0.03cm} \underset{ \vartheta, \mu }{\mathmin} & ~ -\rho(\vartheta) \\
        \mathst~&~ -\mathcal{H}(\xi_h^*; \vartheta, \mu) \leqslant -\epsilon_f, \label{eq:sec4:bilevel-synthesis-outer-constraint-1} \\
        ~&~ -\pi_i(\xi_{\underline{\pi}_i}^*; \mu_i) \leqslant -u_i^{lb} - \epsilon_f , \label{eq:sec4:bilevel-synthesis-outer-constraint-3} \\
        ~&~  \pi_i(\xi_{\overline{\pi}_i}^*; \mu_i) \leqslant u^{ub}_i - \epsilon_f ,\label{eq:sec4:bilevel-synthesis-outer-constraint-4}\\
        ~&~ -s(\xi_s^*) \leqslant -\epsilon_f, \label{eq:sec4:bilevel-synthesis-outer-constraint-2} \\
        ~&~ (\vartheta, \mu) \in \mathcal{O} \defeq \Theta \times M \subset \mathbb{R}^{n_{\vartheta} + n_{\mu}}, \\
        ~&~ \bigl(\xi_h^*,\xi_{\underline{\pi}}^*,\xi_{\overline{\pi}}^*,\xi_s^* \bigr) \in \underset{\xi_h,\xi_{\underline{\pi}},\xi_{\overline{\pi}}, \xi_s}{\argmin} ~ \mathcal{F}(\eta),  \\
        &\hspace{1.2cm} \mathst ~~ -h(\xi_h;\vartheta) \leqslant 0, \label{eq:sec4:bilevel-synthesis-inner-constraint-1}\\
        &\hspace{1.855cm} -h(\xi_{\underline{\pi}_i};\vartheta) \leqslant 0, \label{eq:sec4:bilevel-synthesis-inner-constraint-3}\\
        &\hspace{1.855cm} -h(\xi_{\overline{\pi}_i};\vartheta) \leqslant 0, \label{eq:sec4:bilevel-synthesis-inner-constraint-4}\\
        &\hspace{1.855cm} -h(\xi_s;\vartheta) \leqslant 0, \label{eq:sec4:bilevel-synthesis-inner-constraint-2}\\
        &\hspace{1.855cm} \xi_h,\xi_{\underline{\pi}_i},\xi_{\overline{\pi}_i}, \xi_s \in \mathcal{I} \defeq \mathbb{X} \subset \mathbb{R}^{n}, \label{eq:sec4:bilevel-synthesis-inner-constraint-5}
	\end{align}
\end{subequations}
where $i \in \{1,2, \hdots, m\}$, $\rho:\mathbb{R}^{n_{\vartheta}} \rightarrow \mathbb{R}$ is a mapping function from the parameter vector $\vartheta$ to the size $\rho(\vartheta)$ of $\mathcal{C}(\vartheta)$, and $\epsilon_f \in \mathbb{R}_{>0}$ is a given convergence tolerance for the inner optimization problem (see Theorem \ref{theorem:sec4:optimal-solution-DTCBF} below). The unknown coefficients $\vartheta$ and $\mu$ are selected from compact sets $\Theta \subset \mathbb{R}^{n_\vartheta}$ and $M \subset \mathbb{R}^{n_\mu}$, respectively.
\begin{remark}
    If computing $\rho$ for a parameterized DTCBF is challenging, one may instead optimize an alternative objective (e.g., the Chebyshev radius) that indirectly enlarges $\mathcal{C}(\vartheta)$.
\end{remark}
\begin{remark}
    The inner problem in \eqref{eq:sec4:bilevel-synthesis} may have multiple global minimizers for given outer decision variables $\vartheta$ and $\mu$, making the bilevel problem ill-posed in this case. To address this issue, it is common to categorize bilevel problems as either optimistic or pessimistic \cite{dempe2002foundations}. However, in \eqref{eq:sec4:bilevel-synthesis}, it is evident that the choice of global minimizer in the inner problem does not change the outcome. This is because the outer objective function depends solely on the outer decision variables, and the satisfaction of the outer constraints is determined only by the global minimum of the inner problem. 
\end{remark}

Solving \eqref{eq:sec4:bilevel-synthesis} to global optimality, especially when the inner problem is non-convex for fixed outer decision variables $\vartheta$ and $\mu$, is challenging. Currently, there are two common methods for solving it and achieving $(\epsilon_f, \epsilon_F)$-optimality, given inner optimality tolerance $\epsilon_f \in \mathbb{R}_{> 0}$ and outer optimality tolerance $\epsilon_F \in \mathbb{R}_{> 0}$ (see Definition \ref{def:epsilon_optimality} below). These methods are discussed in Mitsos et al. \cite{Mitsos2008} and Kleniati et al. \cite{Kleniati2014a}, \cite{paulavivcius2020}, with the latter referred to as the Branch-and-Sandwich (B\&S) algorithm. In this paper, 
we use the B\&S algorithm.

\begin{definition}[($\epsilon_f, \epsilon_F$)-optimality \cite{Kleniati2014a}] \label{def:epsilon_optimality}
    Let $\epsilon_f, \epsilon_F \in \mathbb{R}_{> 0}$ be given. Then, $\bar{\eta} \defeq \bigl[\bar{\xi}_h~ \bar{\xi}_{\underline{\pi}}~ \bar{\xi}_{\overline{\pi}}~\bar{\xi}_s~ \bar{\vartheta}~ \bar{\mu}\bigr]^{\T}$ is called an $(\epsilon_f, \epsilon_F)$-optimal solution of the bilevel problem \eqref{eq:sec4:bilevel-synthesis} if $\bar{\eta}$ satisfies the inner and outer constraints and the following inequalities hold%
    \begin{subequations}
    \begin{align}
      \mathcal{H}(\bar{\xi}_h; \bar{\vartheta}, \bar{\mu})
        &\leqslant \mathcal{H}(\xi_h^*; \bar{\vartheta}, \bar{\mu}) + \epsilon_f,
        \label{eq:sec4:prop:valid-DTCBF-bilevel-break-1} \displaybreak[1] \displaybreak[1] \\
      \pi_i(\bar{\xi}_{\underline{\pi}_i}; \bar{\mu}_i)
        &\leqslant \pi_i(\xi_{\underline{\pi}_i}^*; \bar{\mu}_i) + \epsilon_f,
        \label{eq:sec4:prop:valid-DTCBF-bilevel-break-3} \displaybreak[1] \displaybreak[1] \\
      - \pi_i(\bar{\xi}_{\overline{\pi}_i}; \bar{\mu}_i)
        &\leqslant - \pi_i(\xi_{\overline{\pi}_i}^*; \bar{\mu}_i) + \epsilon_f,
        \label{eq:sec4:prop:valid-DTCBF-bilevel-break-4} \displaybreak[1] \displaybreak[1] \\
      s(\bar{\xi}_s) &\leqslant s(\xi_s^*) + \epsilon_f,
        \label{eq:sec4:prop:valid-DTCBF-bilevel-break-2} \displaybreak[1] \displaybreak[1] \\
      -\rho^* &\leqslant \rho(\bar{\vartheta}) + \epsilon_F,
        \label{eq:sec4:prop:valid-DTCBF-bilevel-break-5}
    \end{align}
    \end{subequations}
        for all $i\in \{1,2, \hdots, m\}$, where $[\xi_h^* ~ \xi_{\underline{\pi}}^* ~ \xi_{\overline{\pi}}^* ~ \xi_s^*]^{\T} \in \mathcal{I}^{2+2m}$ corresponds to a global minimizer of the inner problem in \eqref{eq:sec4:bilevel-synthesis} for fixed outer variables $\bar{\vartheta}$ and $\bar{\mu}$, and $\rho^*$ is the global minimum of \eqref{eq:sec4:bilevel-synthesis}.
\end{definition}
\begin{remark}
In the steps of the B\&S algorithm \cite{Kleniati2014a, paulavivcius2020}, we consider the optimization problems \eqref{eq:sec4:proposed-alg:verification-DTCBF}, \eqref{eq:sec4:proposed-alg:admissible-control-input-1}, \eqref{eq:sec4:proposed-alg:admissible-control-input-2}, and \eqref{eq:sec4:proposed-alg:subset-safe-set} individually. This allows us to achieve $\epsilon_f$-optimality, given $\epsilon_f \in \mathbb{R}_{>0}$, for each of them, as in \eqref{eq:sec4:prop:valid-DTCBF-bilevel-break-1}--\eqref{eq:sec4:prop:valid-DTCBF-bilevel-break-2}.
\end{remark}

\begin{theorem} \label{theorem:sec4:optimal-solution-DTCBF}  
Let $\epsilon_f, \epsilon_F \in \mathbb{R}_{>0}$ be given. The function $h$, parameterized by $\bar{\vartheta}$ as in \eqref{eq:sec4:parameterized-DTCBF}, with a given $\gamma \in \mathcal{K}^{\leqslant \mathrm{id}}_{\infty}$, and the control policy $\pi$, parameterized by $\bar{\mu}$ as in \eqref{eq:sec4:parameterized-control-policy}, meet 
Properties \ref{condition:synthesis-1}\textendash\ref{condition:synthesis-3} in Problem \ref{problem:synthesis} (and thus $(h,\gamma)$ with its friend~$\pi$ is a DTCBF) if the coefficient vectors $\bar{\vartheta}$ and $\bar{\mu}$ correspond to an $(\epsilon_f, \epsilon_F)$-optimal solution of the bilevel problem \eqref{eq:sec4:bilevel-synthesis}. Moreover, the zero-superlevel set of $h$ is non-empty, and its size is maximized (within the class of the parameterized DTCBF and the parameterized control policy) up to an $\epsilon_F$ tolerance, i.e., $\rho(\bar{\vartheta}) + \epsilon_F \geqslant -\rho^*$, where $\rho^*$ is the global minimum of \eqref{eq:sec4:bilevel-synthesis}.

    \begin{proof}
        Assume $\bar{\eta} \defeq \bigl[\bar{\xi}_h~ \bar{\xi}_{\underline{\pi}}~ \bar{\xi}_{\overline{\pi}}~\bar{\xi}_s~ \bar{\vartheta}~ \bar{\mu}\bigr]^{\T}$ is an $(\epsilon_f, \epsilon_F)$-optimal solution of \eqref{eq:sec4:bilevel-synthesis} for given $\epsilon_f, \epsilon_F \in \mathbb{R}_{>0}$. Based on \mbox{Definition \ref{def:epsilon_optimality}}, $\bar{\eta}$ satisfies the inner and outer constraints of the bilevel problem \eqref{eq:sec4:bilevel-synthesis} and also satisfies \eqref{eq:sec4:prop:valid-DTCBF-bilevel-break-1}\textendash\eqref{eq:sec4:prop:valid-DTCBF-bilevel-break-5}. The satisfaction of the outer constraint \eqref{eq:sec4:bilevel-synthesis-outer-constraint-1} and \eqref{eq:sec4:prop:valid-DTCBF-bilevel-break-1} implies that $0 \leqslant \mathcal{H}(x; \bar{\vartheta}, \bar{\mu})$ for all $x \in \mathcal{C}(\bar{\vartheta})$; that is, Property~\ref{condition:synthesis-1} is satisfied. Similarly, \eqref{eq:sec4:bilevel-synthesis-outer-constraint-3} and \eqref{eq:sec4:bilevel-synthesis-outer-constraint-4}, along with \eqref{eq:sec4:prop:valid-DTCBF-bilevel-break-3} and \eqref{eq:sec4:prop:valid-DTCBF-bilevel-break-4}, for all $i \in \{1,2, \hdots, m\}$, imply that
        $u^{lb}_i \leqslant \pi_i(x; \bar{\mu}_i) \leqslant u^{ub}_i$ holds for all $i \in \{1,2, \hdots, m\}$ and for all $x \in \mathcal{C}(\bar{\vartheta})$; that is Property~\ref{condition:synthesis-2} is met. Moreover, \eqref{eq:sec4:bilevel-synthesis-outer-constraint-2} and \eqref{eq:sec4:prop:valid-DTCBF-bilevel-break-2} imply that $0 \leqslant s(x)$ holds for all $x \in \mathcal{C}(\bar{\vartheta})$; that is, Property~\ref{condition:synthesis-3} is met. Thus, the function $h$, parameterized by $\bar{\vartheta}$ as in \eqref{eq:sec4:parameterized-DTCBF}, with given $\gamma$ and its friend $\pi$, parameterized by $\bar{\mu}$ as in \eqref{eq:sec4:parameterized-control-policy}, is a DTCBF.
        
        Additionally, satisfying \eqref{eq:sec4:prop:valid-DTCBF-bilevel-break-5} and the inner constraints in \eqref{eq:sec4:bilevel-synthesis} implies that $\mathcal{C}(\bar{\vartheta})$ is non-empty and its size is maximized within an $\epsilon_F \in \mathbb{R}_{>0}$ tolerance.
    \end{proof}
\end{theorem}
\begin{remark}
The bilevel problem \eqref{eq:sec4:bilevel-synthesis} is feasible if there exist coefficient vectors $(\bar{\vartheta},\bar{\mu}) \in \mathcal{O} \defeq \Theta \times M$ such that the function $h$, parameterized by $\bar{\vartheta}$ as in \eqref{eq:sec4:parameterized-DTCBF}, with a given $\gamma \in \mathcal{K}_{\infty}^{\leqslant \mathrm{id}}$, and the control policy $\pi$, parameterized by $\bar{\mu}$ as in \eqref{eq:sec4:parameterized-control-policy}, satisfy 
\begin{align}
    \epsilon_f & \leqslant \mathcal{H}(x; \bar{\vartheta}, \bar{\mu}), \nonumber \\
    u^{lb}_i + \epsilon_f &\leqslant \pi_i(x; \bar{\mu}_i) \leqslant u^{ub}_i - \epsilon_f, \text{~~for all~} i \in \{1,2, \hdots, m\}, \nonumber \\
    \epsilon_f &\leqslant s(x), \nonumber
\end{align}
for all $x \in \{ \xi \in \mathbb{R}^n \mid h(\xi) \geqslant 0\}$ (that is, Properties \ref{condition:synthesis-1}\textendash\ref{condition:synthesis-3} are satisfied \textit{with conservatism level $\epsilon_f$}). As a result, it is advisable to select $\epsilon_f$ sufficiently small to avoid infeasibility issue and to reduce the conservatism of the synthesized DTCBF $h$ and the synthesized control policy~$\pi$.
\end{remark}
\begin{theorem} \label{theorem:sec4:MFCQ-inner}
    Let $\epsilon_f, \epsilon_F \in \mathbb{R}_{>0}$ be given. Under Assumptions \ref{assump:sec4:bounded-DTCBF}\textendash\ref{assump:sec4:MFCQ}, the B\&S algorithm \cite{Kleniati2014a, paulavivcius2020} computes an $(\epsilon_f, \epsilon_F)$-optimal solution of \eqref{eq:sec4:bilevel-synthesis} in finite time.
    
\begin{proof} 
    Since $\mathcal{C}(\vartheta) \subset \mathbb{X}$ and $\partial \mathcal{C}(\vartheta) \cap \partial \mathbb{X} = \emptyset$ hold for all $\vartheta \in \Theta$ due to Assumption \ref{assump:sec4:chosen-X}, the linear inequality constraints \eqref{eq:sec4:bilevel-synthesis-inner-constraint-5} and the possibly nonlinear constraints \eqref{eq:sec4:bilevel-synthesis-inner-constraint-1}--\eqref{eq:sec4:bilevel-synthesis-inner-constraint-2} cannot be simultaneously active. Thus, because \mbox{$\nabla_x h(x; \vartheta) \neq [0~ \hdots~ 0]^{\T} \in \mathbb{R}^n$} for all $\vartheta \in \Theta$ and for all $x \in \mathbb{X}$ such that $h(x; \vartheta) = 0$ (by Assumption \ref{assump:sec4:MFCQ}), it is straightforward to show that the MFCQ \cite{Bertsekas2016} holds for the inner problem of \eqref{eq:sec4:bilevel-synthesis} at all \mbox{$\eta \in \mathcal{I}^{2+2m} \times \Theta \times M$}. As a result, under Assumptions \ref{assump:sec4:bounded-DTCBF}\textendash\ref{assump:sec4:MFCQ} and given $\epsilon_f, \epsilon_F \in \mathbb{R}_{> 0}$, the B\&S algorithm computes an $(\epsilon_f, \epsilon_F)$-optimal solution of \eqref{eq:sec4:bilevel-synthesis} in a finite number of iterations (see \cite{Kleniati2014a} for more details).
\end{proof}
\end{theorem}
\subsection{Strategies for Simplifying the Bilevel Problem} \label{sec:sec4:strategies-simplifying}
The computation time for solving the bilevel problem \eqref{eq:sec4:bilevel-synthesis} may be high, particularly for high-dimensional dynamical systems. Specifically, adding decision variables in \eqref{eq:sec4:bilevel-synthesis} causes the computation time to increase exponentially. In this subsection, our aim is to introduce strategies to simplify meeting Properties \ref{condition:synthesis-2} and \ref{condition:synthesis-3}, thereby making \eqref{eq:sec4:bilevel-synthesis} more tractable.

\begin{remark}
    Choosing a sufficiently large $\epsilon_F \in \mathbb{R}_{>0}$ (the tolerance for maximizing the zero-superlevel set of the resulting DTCBF) significantly reduces the computation time of the bilevel problem \eqref{eq:sec4:bilevel-synthesis}. Additionally, if finding a valid but potentially conservative DTCBF is sufficient, one can terminate the B\&S algorithm once a feasible solution, one that satisfies the inner and outer constraints of \eqref{eq:sec4:bilevel-synthesis}, is found.
\end{remark}
\subsubsection{Property \ref{condition:synthesis-2}}
In the proposed bilevel problem \eqref{eq:sec4:bilevel-synthesis}, the admissibility of each parameterized control policy $\pi_i$, \mbox{$i \in \{1, \hdots, m\}$}, is verified for all $x \in \mathcal{C}(\vartheta)$ by imposing \eqref{eq:sec4:proposed-alg:admissible-control-input-1-ineq} and \eqref{eq:sec4:proposed-alg:admissible-control-input-2-ineq}. This approach may become computationally expensive, especially for systems with multiple control inputs. To make this approach more tractable, we can approximate $\mathbb{U} \defeq \left[u^{lb}, ~ u^{ub}\right]$ by an $m$-dimensional super-ellipsoid. Specifically, we define an inner approximation of $\mathbb{U}$ as
\begin{align} \label{eq:sec4:strategies:approx-input}
    \widetilde{\mathbb{U}} & \defeq \Biggl\{ u \in \mathbb{R}^m \, \Bigg| \, \sum_{i=1}^{m} \biggl| \frac{u_i - \frac{u_i^{ub} + u_i^{lb}}{2}}{(u_i^{ub} - u_i^{lb})/2} \biggr|^p \leqslant 1 \Biggr\},
\end{align}
where $p \in \mathbb{N}$ and $|.|$ is the absolute value. We can show that $\widetilde{\mathbb{U}} \subset \mathbb{U}$ and $\widetilde{\mathbb{U}}$ converges to $\mathbb{U}$ as $p$ increases. Hence, we assert that the control policy $\pi$, parameterized by $\mu$, is admissible with respect to $\mathbb{U}$ for all $x \in \mathcal{C}(\vartheta)$ if
\begin{align}
\sum_{i=1}^{m} \biggl| \frac{\pi_i(\xi_{\pi}^*;\mu_i) - \frac{u_i^{ub} + u_i^{lb}}{2}}{(u_i^{ub} - u_i^{lb})/2} \biggr|^p \leqslant 1, \nonumber
\end{align}
where $\xi_\pi^* \in \mathbb{R}^n$ is a global minimizer of 
\begin{subequations} \nonumber
\begin{align}
     \xi_\pi^* \in \underset{\xi_\pi \in \mathbb{X}}{\argmin} ~&~ -\sum_{i=1}^{m} \biggl| \frac{\pi_i(\xi_{\pi};\mu_i) - \frac{u_i^{ub} + u_i^{lb}}{2}}{(u_i^{ub} - u_i^{lb})/2} \biggr|^p \\
    \mathst ~&~ -h(\xi_\pi; \vartheta) \leqslant 0.
\end{align}
\end{subequations}
Although this formulation yields a more conservative problem than \eqref{eq:sec4:bilevel-synthesis}, it reduces the number of inner variables associated with verifying the admissibility of the control policy from $2nm$ to $n$, where $n, m \in \mathbb{N}$ denote the number of states and the number of control inputs in the system \eqref{eq:sec2:dynamical-system}, respectively.

If the control admissible set \mbox{$\mathbb{U} \defeq [-\bar{u},~ \bar{u}]$}, with \mbox{$\bar{u} \in \mathbb{R}^m$}, is symmetric about zero, an alternative approach to verifying the admissibility of the control policy is to check whether 
 $\pi_i(\xi^*_{\pi_i}; \mu_i)^2 \leqslant \bar{u}_i^2$
holds for all $i \in \{1, \hdots, m\}$, where \mbox{$\xi^*_{\pi_i} \in \mathbb{R}^n$} is a global minimizer of 
\begin{subequations}
\begin{align}
   \xi^*_{\pi_i} \in \underset{ \xi_{\pi_i} \in \mathbb{X} }{\argmin} ~&~  - \pi_i(\xi_{\pi_i}; \mu_i)^2 \nonumber \\
   \mathst ~&~ -h(\xi_{\pi_i}; \vartheta) \leqslant 0. \nonumber
\end{align}
\end{subequations}
In this approach, without any conservatism, the number of required inner variables for verifying the admissibility of the control policy is reduced from $2nm$ to $nm$.
\subsubsection{Property \ref{condition:synthesis-3}}
Instead of imposing \eqref{eq:sec4:proposed-alg:subset-safe-set-ineq} and obtaining a global minimizer of \eqref{eq:sec4:proposed-alg:subset-safe-set} to ensure that Property \ref{condition:synthesis-3}, $\mathcal{C}(\vartheta) \subseteq \mathcal{S}$, holds, we can directly impose constraints on $\vartheta$. For instance, suppose $h$ is parameterized in such a way that its zero-superlevel set forms an ellipse. In this scenario, we can impose constraints on the lengths of its semi-minor and semi-major axes to ensure the ellipse remains within the safe set $\mathcal{S}$.

\section{Numerical Case Study} 
\label{sec:sec5}
In this section, we will demonstrate our methods through numerical case studies.\footnote{The code for the numerical case studies is available at \href{https://github.com/shakhesierfan/Verification-and-Synthesis-of-DTCBFs}{github.com/shakhesierfan/Verification-and-Synthesis-of-DTCBFs}.}
The subsequent results are obtained on a laptop with an Intel i7-12700H processor.
\subsection{Synthesis} \label{sec:sec5:synthesis}
\subsubsection{Polynomial System} \label{sec:sec5:synthesis-polynomial}
We discretize the continuous-time system in \cite{Wang2023a} using the forward Euler discretization with a sample time of $T_s \defeq \unit[0.1]{s}$, leading to
\begin{align} \label{eq:sec5:example:discrete-dynamical}
    \begin{bmatrix}
        x^+_{1} \\ x^+_{2}
    \end{bmatrix} = \begin{bmatrix}
        x_1 + x_2T_s + (x_1^2 + x_2 + 1)T_su_1 \\
        x_2 + (x_1 + \frac{1}{3}x_1^3 + x_2)T_s + (x_2^2 + x_1 + 1)T_su_2  
    \end{bmatrix},
\end{align}
where $u_1 \in [-1.5,~ 1.5]$ and $u_2 \in [-1.5,~ 1.5]$. The discrete-time system \eqref{eq:sec5:example:discrete-dynamical} can concisely be written as $x^+ = f(x, u)$. The safe set $\mathcal{S}$ is defined as the disk 
\begin{align}  \label{eq:sec5:example:safe-set}
    \mathcal{S} \defeq \{x \in \mathbb{R}^2 \mid x_1^2 + x_2^2 - 3 \leqslant 0 \}.
\end{align} 
In this example, we consider $\gamma \in \mathcal{K}^{\leqslant \mathrm{id}}_{\infty}$ as a linear function, namely $\gamma(r) \defeq 0.5r$, $r \in \mathbb{R}_{\geqslant 0}$.

We aim to synthesize a function $h$ and a control policy $\pi$ for the system \eqref{eq:sec5:example:discrete-dynamical} and the safe set \eqref{eq:sec5:example:safe-set}, such that $(h, \gamma)$ with its friend $\pi$ satisfies Properties \ref{condition:synthesis-1}\textendash\ref{condition:synthesis-3} in Problem \ref{problem:synthesis}, consequently, $(h, \gamma)$ with its friend $\pi$ is a DTCBF. To this end, we parameterize the function $h$, such that its zero-superlevel set $\mathcal{C}(\vartheta)$ is ellipsoidal and centered at the origin, as
\begin{align} \label{eq:sec5:example:parameterized-DTCBF}
    h(x; \vartheta) \defeq 1 -\vartheta_1x_1^2 - \vartheta_2x_1x_2 - \vartheta_3x_2^2,
\end{align}
where $\vartheta \defeq [\vartheta_1~ \vartheta_2~ \vartheta_3]^{\T} \in \mathbb{R}^3$ is unknown. The ellipse $h(x; \vartheta) = 0$ is a non-degenerate real ellipse if and only if $L(\vartheta) \defeq \vartheta_2^2 - 4\vartheta_1\vartheta_3 < 0$, $\vartheta_1 > 0$, and $\vartheta_3 > 0$. Additionally, the length of the semi-major axis $a(\vartheta)$ and the length of the semi-minor axis $b(\vartheta)$ of the ellipse are computed as
\begin{align}
    a(\vartheta), b(\vartheta)  = \frac{-\sqrt{-2L(\vartheta)\left(\vartheta_1 + \vartheta_3 \pm \sqrt{(\vartheta_1 - \vartheta_3)^2 + \vartheta_2^2}\right)}}{L(\vartheta)}. \nonumber
\end{align} 
Thus, we can assert that the zero-superlevel set $\mathcal{C}(\vartheta)$ of $h$ is a subset of the safe set, $\mathcal{C}(\vartheta) \subseteq \mathcal{S}$, if 
    \mbox{$a(\vartheta)^2 \leqslant 3$} and \mbox{$b(\vartheta)^2 \leqslant 3$}.
Moreover, the area $\rho(\vartheta)$ of $\mathcal{C}(\vartheta)$ is calculated as 
\begin{align}
    \rho(\vartheta) = \pi a(\vartheta)b(\vartheta) = \frac{2\pi}{\sqrt{-L(\vartheta)}}. \nonumber
\end{align}
Therefore, to maximize the size of $\mathcal{C}(\vartheta)$, it suffices to maximize the value of $L(\vartheta)$.

\begin{figure}[t]\centering 
\includegraphics[width=0.42\textwidth]{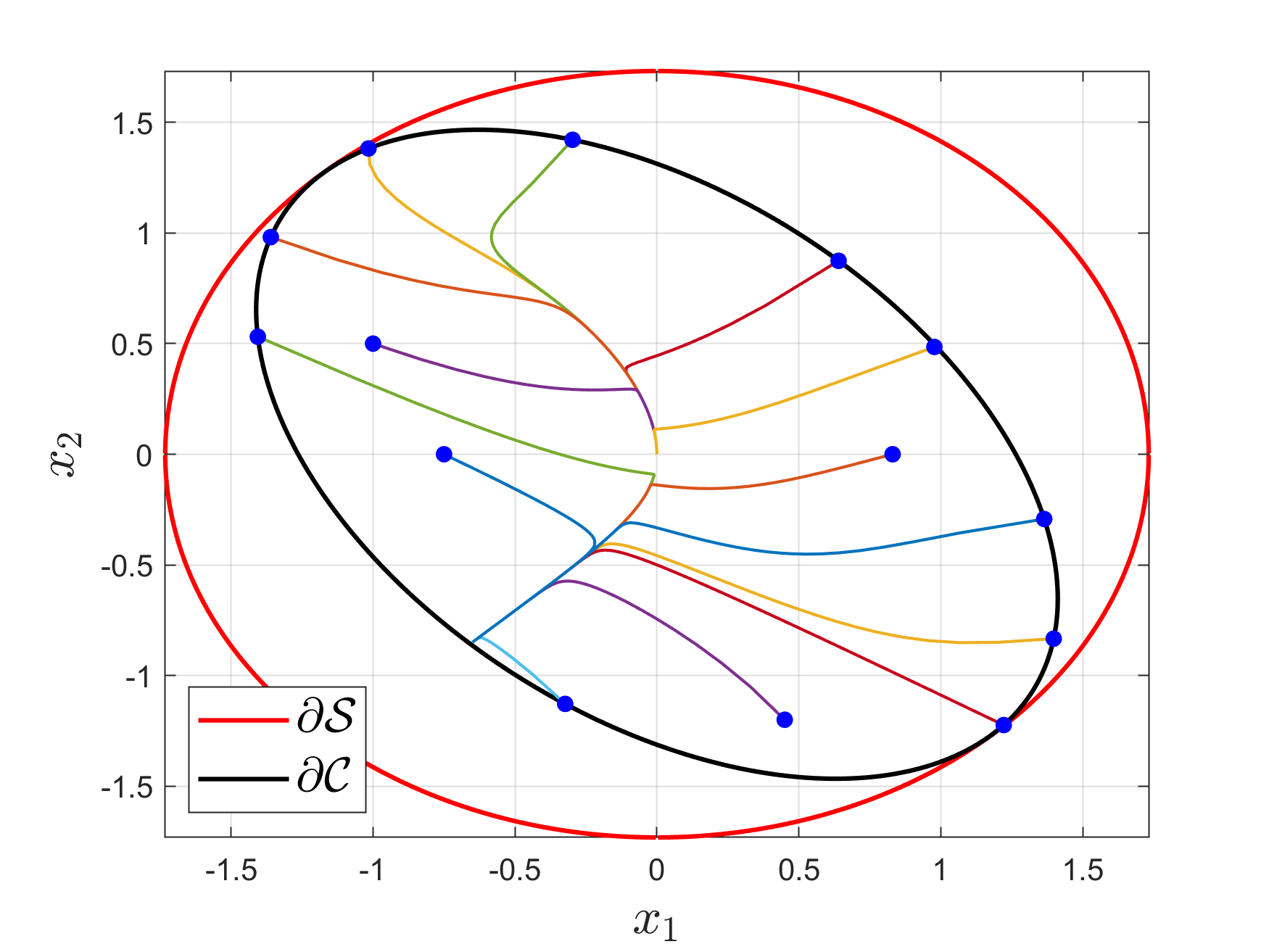}
    \caption{The black ellipse represents the boundary of the zero-superlevel set of the synthesized DTCBF for the system \eqref{eq:sec5:example:discrete-dynamical} and the safe set \eqref{eq:sec5:example:safe-set} (red circle). The colored curves represent trajectories that originate from the blue filled circles and evolve according to the dynamics when applying the synthesized control policy.}
    \label{fig:sec4:synthesis-example}
\end{figure}

We also parameterize the control policy $\pi$ as 
\begin{subequations} \label{eq:sec5:example:control-policy}
\begin{align}
    \pi_1(x; \mu_1) &\defeq \mu_{1,1}x_1 + \mu_{1,2}x_2, \label{eq:sec5:example:control-policy-1} \\
    \pi_2(x; \mu_2) &\defeq \mu_{2,1}x_1 + \mu_{2,2}x_2, \label{eq:sec5:example:control-policy-2}
\end{align}
\end{subequations}
with $\mu_1 \defeq [\mu_{1,1} ~ \mu_{1,2}]^{\T} \in \mathbb{R}^2$, $\mu_2 \defeq [\mu_{2,1} ~ \mu_{2,2}]^{\T} \in \mathbb{R}^2$, $\mu \defeq [\mu_1 ~ \mu_2]^{\T}$, and $\pi(x; \mu) \defeq [\pi_1(x; \mu_1) ~\pi_2(x; \mu_2)]^{\T}$.
Since the control admissible set is symmetric about zero, the control policy $\pi$ in \eqref{eq:sec5:example:control-policy} is admissible for all $x \in \mathcal{C}(\vartheta)$ if $\pi_1(\xi_{\pi_1}^*;\mu_1)^2 \leqslant 1.5^2$ and $\pi_2(\xi_{\pi_2}^*;\mu_2)^2 \leqslant 1.5^2$,
where \mbox{$\xi_{\pi_i}^* \in \mathbb{R}^2$}, $i \in \{1, 2\}$, is a global minimizer of
\begin{subequations}  \nonumber
    \begin{align}
     \xi_{\pi_i}^* \in \underset{ \xi_{\pi_i} \in \mathbb{X} }{\argmin} ~&~ -\pi_i(\xi_{\pi_i};\mu_i)^2 \\
    \mathst ~&~ -h(\xi_{\pi_i}; \vartheta) \leqslant 0.
    \end{align}
\end{subequations} 
In this example, we select the inner variable domain $\mathcal{I} \defeq \mathbb{X}$ such that $\mathcal{S}\subset \mathcal{I}$ as 
\begin{align}
    \mathcal{I} \defeq \mathbb{X} \defeq [-1.75,~ 1.75] \times [-1.75,~ 1.75], \nonumber 
\end{align}
and the outer variable domain as $\mathcal{O} \defeq \Theta \times M$ where 
\begin{align}
\vartheta &\in \Theta \hspace{0.1cm} \defeq [0.1, ~ 1.5]\times [-1.5, ~ 1.5] \times [0.1, ~ 1.5], \nonumber \\
\mu &\in M \defeq [-5, ~ -0.1]^4. \nonumber
\end{align} 
As a result, the bilevel optimization problem to synthesize $h$ for this example is constructed as 
\begin{subequations} \label{eq:sec5:example-bilevel}
	\begin{align}
		\underset{ \vartheta, \mu }{\mathmin} ~&~ -L(\vartheta) \\
        \mathst ~&~ -\mathcal{H}(\xi_h^*; \vartheta, \mu) \leqslant -\epsilon_f,  \\
        ~&~ \pi_1(\xi_{\pi_1}^*; \mu_1)^2 \leqslant 1.5^2 - \epsilon_f,  \\
        ~&~ \pi_2(\xi_{\pi_2}^*; \mu_2)^2 \leqslant 1.5^2 - \epsilon_f,  \\
        ~&~ a(\vartheta)^2 \leqslant 3, \\
        ~&~ b(\vartheta)^2 \leqslant 3, \\
        ~&~ L(\vartheta) \leqslant -10^{-3}, \\
        ~&~ (\vartheta, \mu) \in \Theta \times M, \\
        ~&~ \left(\xi_h^*, \xi_{\pi_1}^*, \xi_{\pi_2}^*\right) \in \underset{\xi_h, \xi_{\pi_1}, \xi_{\pi_2}}{\argmin} ~ \mathcal{F}(\xi_h, \xi_{\pi_1}, \xi_{\pi_2}; \vartheta, \mu),  \\
        &\hspace{1.2cm} \mathst ~ -h(\xi_h;\vartheta) \leqslant 0, \\
        &\hspace{1.745cm} -h(\xi_{\pi_1};\vartheta) \leqslant 0, \\
        &\hspace{1.745cm} -h(\xi_{\pi_2};\vartheta) \leqslant 0, \\
        &\hspace{1.745cm} \xi_h,\xi_{\pi_1}, \xi_{\pi_2}  \in \mathcal{I},
	\end{align}
\end{subequations}
where we define 
\begin{align}
    \mathcal{F}(\xi_h, \xi_{\pi_1}, \xi_{\pi_2}; \vartheta, \mu) &\defeq \mathcal{H}(\xi_h; \vartheta, \mu) \nonumber \\
    &\qquad - \pi_1(\xi_{\pi_1}; \mu_1)^2 - \pi_2(\xi_{\pi_2}; \mu_2)^2, \nonumber
\end{align}
and
\begin{align}
    \mathcal{H}(\xi_h; \vartheta, \mu) \defeq h(f(\xi_h, \pi(\xi_h; \mu)); \vartheta) - (1-0.5)h(\xi_h; \vartheta). \nonumber
\end{align}
For \eqref{eq:sec5:example-bilevel}, it can be shown that \mbox{Assumptions \ref{assump:sec4:bounded-DTCBF}\textendash\ref{assump:sec4:MFCQ}} are met. Thus, we apply the B\&S algorithm to obtain \mbox{$(\epsilon_f, \epsilon_F)$-optimality} of \eqref{eq:sec5:example-bilevel} for the inner optimality tolerance $\epsilon_f \defeq 0.001$ and the outer optimality tolerance \mbox{$\epsilon_F \defeq 0.4$}. The sub-problems of the B\&S algorithm are solved using the BARON solver \cite{Sahinidis1996}. After around 1000 iterations, where each iteration takes an average of \unit[4]{s} (total time about $\unit[66]{min}$), an $(\epsilon_f, \epsilon_F)$-optimal solution of \eqref{eq:sec5:example-bilevel} is obtained as 
\begin{align}
    \vartheta_1 &= 0.626, \hspace{0.6cm} \vartheta_2 = 0.537, \hspace{0.38cm}\vartheta_3 = 0.580, \nonumber \\
    \mu_{1,1} &= -0.976, ~ \mu_{1,2} = -1, \nonumber \\ \mu_{2,1} &= -0.976, ~ \mu_{2,2} = -1. \nonumber
\end{align}
It should be noted that this solution is already obtained after approximately 150 iterations (about $\unit[10]{min}$). However, to ensure that there is no DTCBF with a larger zero-superlevel set within $\epsilon_F \defeq 0.4$ tolerance, approximately 850 more iterations are required in the B\&S algorithm to discard other nodes. The synthesized DTCBF is depicted in Figure \ref{fig:sec4:synthesis-example}.
\subsubsection{Cart-Pole} \label{sec:sec5:synthesis-cart-pole}
Consider the discretized cart-pole system discussed in \cite{Neural-Lyapunov} as 
\begin{align} \label{eq:dynamics-cart-pole}
    \begin{bmatrix}
        x_c^+ \\
        v_c^+ \\
        \theta^+ \\
        \omega^+
    \end{bmatrix} \hspace{-0.1cm}= \hspace{-0.1cm}\begin{bmatrix}
        x_c + v_cT_s \\
        v_c + \frac{T_s\left(u + m_p \sin\theta \hspace{0.02cm}(L\omega^2 - g\cos \theta)\right)}{m_c + m_p\sin^2 \theta} \\
        \theta + \omega T_s \\
        \omega + \frac{T_s\left(-u\cos \theta - m_pL\omega^2 \cos\theta \hspace{0.02cm} \sin \theta + (m_c + m_p)g \sin \theta\right)}{L(m_c + m_p\sin^2 \theta)}
    \end{bmatrix},
\end{align}
where $x_c, v_c \in \mathbb{R}$ are the horizontal position and velocity of the cart, respectively, $\theta \in \mathbb{R}$ is the angle of the pole from the upward vertical direction, and $\omega \in \mathbb{R}$ is its angular velocity. Moreover, $x \defeq [x_c ~ v_c ~ \theta ~ \omega]^{\T}$, and \mbox{$u \in \mathbb{U} \defeq [-25, ~ 25]$} is the horizontal force applied to the cart, which is considered as the control input to the system. We consider the following parameters: the cart mass $m_c \defeq \unit[2]{kg}$, the pole mass \mbox{$m_p \defeq \unit[0.1]{kg}$}, the pole length $L \defeq \unit[1]{m}$, and the sample time $T_s \defeq \unit[0.01]{s}$.
The safe set $\mathcal{S}$, shown in Figure \ref{fig:sec4:verif-cart-pole-example}, is defined as
\begin{align} \label{eq:safe-set-cart-pole}
    \mathcal{S} \defeq \{ x \in \mathbb{R}^4 \bigm|\theta^2 + \omega^2 \leqslant (\pi/4)^2 \}.
\end{align}
We aim to synthesize a function $h$ and a control policy $\pi$ for the system \eqref{eq:dynamics-cart-pole} and the safe set \eqref{eq:safe-set-cart-pole}, such that $(h, \gamma)$ with its friend $\pi$ is a DTCBF. We parameterize $h$ as
\begin{align}
    h(x; \vartheta) = 1 + \vartheta_1\omega^2 + \vartheta_2 \theta^2 + \vartheta_3 \omega \theta + \vartheta_4 \omega + \vartheta_5 \theta,
\end{align}
where $\vartheta = [\vartheta_1~ \hdots~ \vartheta_5]^{\T} \in \Theta \defeq [-5,~ -0.1]^2\times [-2,~ -0.1] \times [0.1,~ 1]^2$, and $\pi$ as
\begin{align}
    \pi(x; \mu) = \mu_1\omega^2 + \mu_2 \theta^2 + \mu_3 \omega \theta + \mu_4 \omega + \mu_5 \theta,
\end{align}
where $\mu = [\mu_1 ~\hdots ~\mu_5]^{\T} \in M \defeq [-30,~ 30]^5$. We consider $\gamma \in \mathcal{K}_{\infty}^{\leqslant \mathrm{id}}$ as the identity function and define the inner variable domain as $\mathcal{I} \defeq [-\pi/4,~ \pi/4]^2$. 
We then apply the bilevel approach proposed in Section~\ref{sec:sec4}, with \mbox{$\epsilon_f \defeq 10^{-5}$} and \mbox{$\epsilon_F \defeq 0.5$}, without repeating its formulation here. The sub-problems of the B\&S algorithm are solved using the Couenne solver \cite{Belotti01102009}. After around 600 iterations (total time about $\unit[80]{min}$), the DTCBF $h$, shown in Figure~\ref{fig:sec4:verif-cart-pole-example}, and the control policy $\pi$ are synthesized as 
\begin{align}
    h(x) &= -2.5\omega^2 - 2.5\theta^2 -1.7 \omega \theta + 0.1 \omega + 0.1\theta + 1, \label{sec5-2:synthesized-DTCBF} \\
    \pi(x) &= 16.2 \omega^2 -4.8 \theta^2 + 18.0 \omega \theta + 16.4 \omega + 30.0 \theta.  \label{sec5-2:synthesized-control-policy} 
\end{align}

\subsection{Verification} \label{sec:sec5:verification}
\subsubsection{Linear System}
\begin{figure*}[ht!] 
\begin{subfigure}{.33\textwidth}
\centering
\includegraphics[width=1\textwidth]{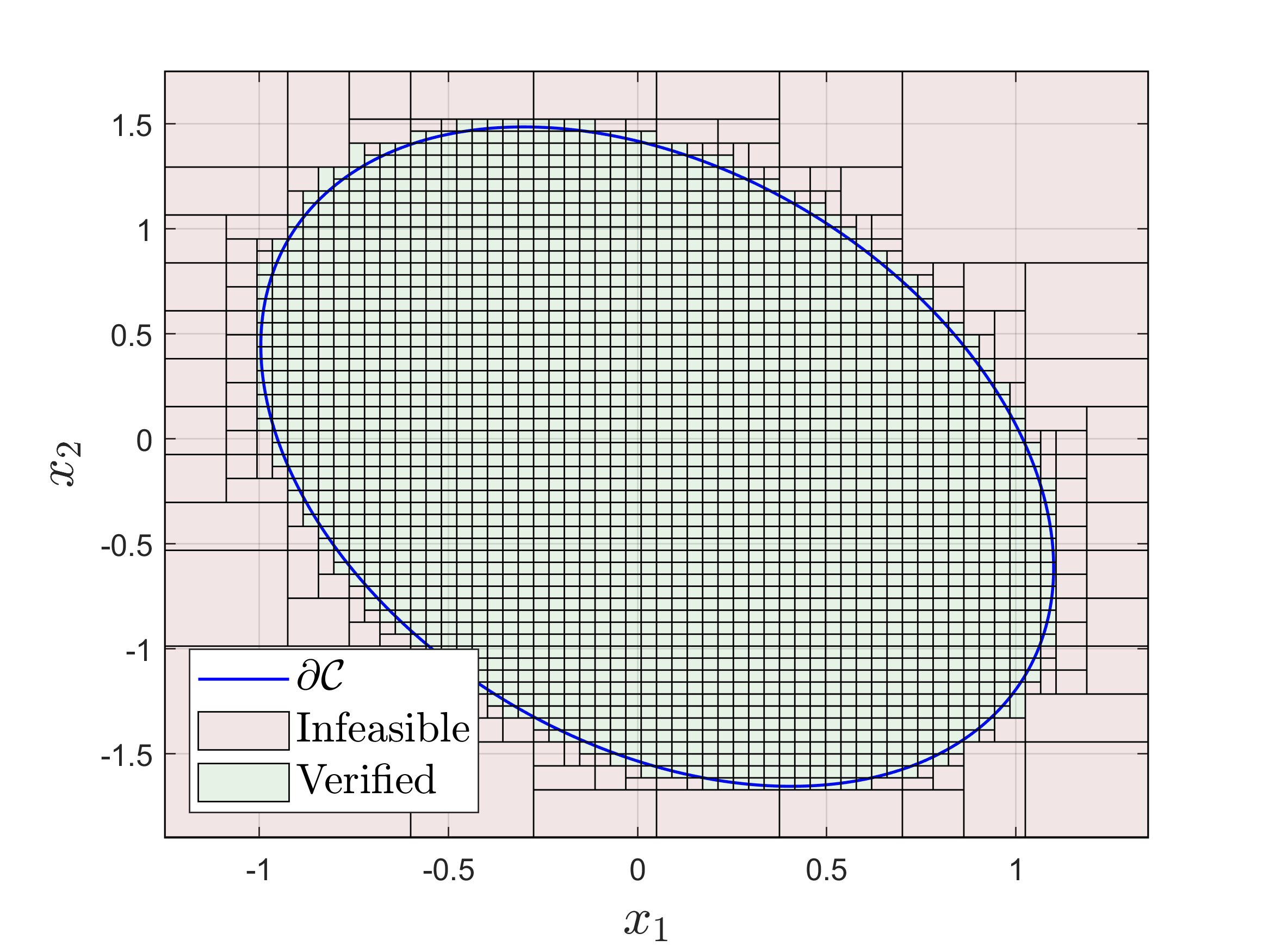} \caption{} \label{fig:example-DTCBF:1a}
\end{subfigure}
\begin{subfigure}{.33\textwidth}
\centering
\includegraphics[width=1\textwidth]{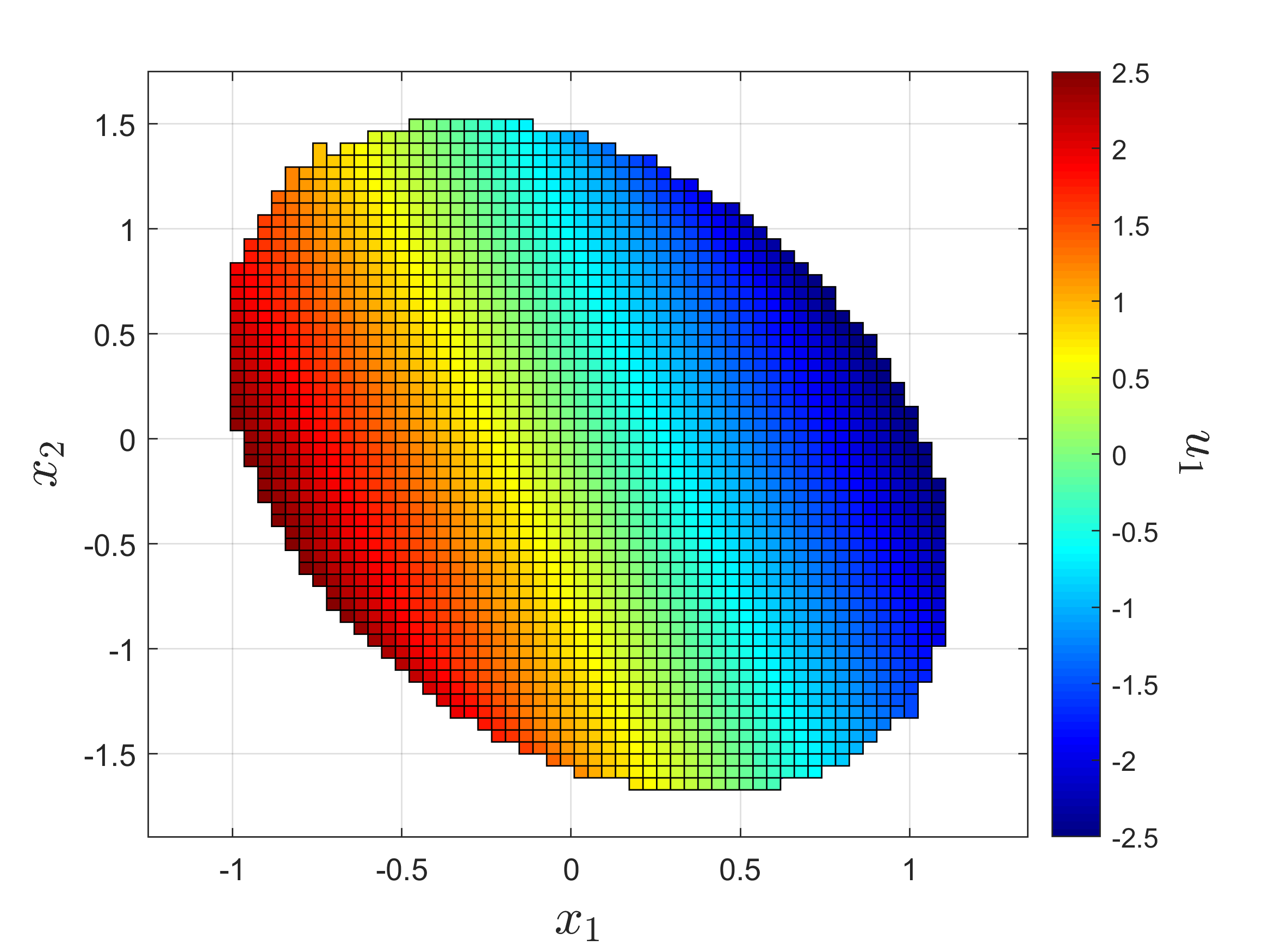} \caption{} \label{fig:example-DTCBF:1b}
\end{subfigure}
\begin{subfigure}{.33\textwidth}
\centering
\includegraphics[width=1\textwidth]{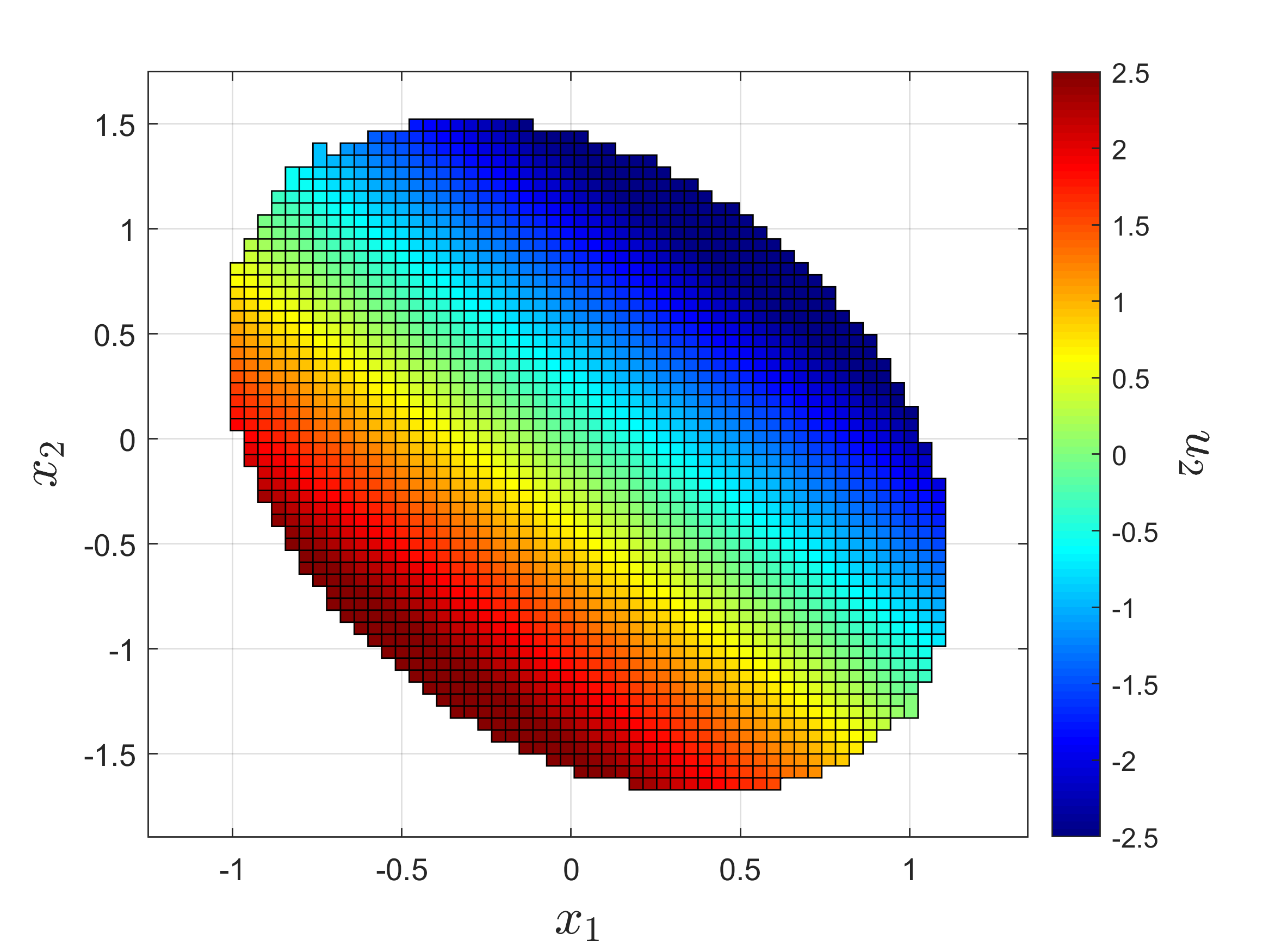} \caption{} \label{fig:example-DTCBF:1c}
\end{subfigure}
\caption{Algorithm \ref{algorithm:verification} is applied to the candidate DTCBF \eqref{eq:sec5:example:DTCBF} with $\gamma$ defined as $\gamma(r) \defeq 0.8r$, $r \in \mathbb{R}_{\geqslant 0}$, for the discrete-time system \eqref{eq:sec5:example:verification:discrete-dynamical}. \mbox{(a) The green rectangles} signify that $(h,\gamma)$ is verified within the corresponding subdomain, Case \ref{case:convex-verification-known-A}, while the red rectangles indicate that the optimization problem \eqref{eq:sec3:verification:unknown:step3-convex} is infeasible, Case \ref{case:convex-verification-known-B}, in the sense that the corresponding subdomain is entirely outside $\mathcal{C}$. The blue curve represents the boundary of $\mathcal{C}$. (b, c) The colored subdomains represent the computed values of the control inputs on each subdomain, namely, $u_1$ in (b) and $u_2$ in (c).}
\label{fig:sec5:verification-example}
\end{figure*}
We consider the continuous-time dynamical system discussed in \cite{wang2022}:
\begin{align} \label{eq:sec5:example:verification:continuous-dynamical}
    \begin{bmatrix}
        \dot{x}_1 \\ \dot{x}_2
    \end{bmatrix} = \begin{bmatrix}
        2 & 1 \\
        3 & 1
    \end{bmatrix}\begin{bmatrix}
        x_1 \\ x_2
    \end{bmatrix} + \begin{bmatrix}
        u_1 \\ u_2
    \end{bmatrix},
\end{align}
where $u_1 \in [-2.5,~ 2.5]$ and $u_2 \in [-2.5,~ 2.5]$. In \cite{wang2022}, the continuous-time CBF is synthesized as 
\begin{align} \label{eq:sec5:example:DTCBF}
    h(x) &= −7.635x_1^2 - 3.439x_1x_2 − 3.4024x_2^2 \nonumber \\
    &\quad + 0.5x_1 − 0.4x_2 + 7.402,
\end{align}
and its corresponding control policy as 
\begin{subequations} \label{eq:sec5:example:Control-Policy}
\begin{align}
    \pi_1(x) &= −2.32x_1 − 1.11x_2 + 0.022,\\
    \pi_2(x) & = −2.12x_1 − 1.27x_2 − 0.046,
\end{align}
\end{subequations}
with $\pi(x) \defeq [\pi_1(x)~ \pi_2(x)]^{\T}$.
Consider the discrete-time representation of \eqref{eq:sec5:example:verification:continuous-dynamical} using exact zero-order hold (ZOH) discretization with a sample time of $T_s \defeq \unit[1]{s}$, leading to 
\begin{align} \label{eq:sec5:example:verification:discrete-dynamical}
    \begin{bmatrix}
        x^+_{1} \\ x^+_{2}
    \end{bmatrix} = \begin{bmatrix}
        17.6 & 7.3 \\
        22.0 & 10.3
    \end{bmatrix}\begin{bmatrix}
        x_{1} \\ x_{2}
    \end{bmatrix} + \begin{bmatrix}
        5.4 & 2.0 \\
        5.9 & 3.4
    \end{bmatrix} \begin{bmatrix}
        u_{1} \\ u_{2}
    \end{bmatrix},
\end{align}
where $u_1 \in [-2.5,~ 2.5]$ and $u_2 \in [-2.5,~ 2.5]$. We aim to verify whether $(h, \gamma)$, where $\gamma(r) \defeq 0.8r$, \mbox{$r \in \mathbb{R}_{\geqslant 0}$}, with the control policy $\pi$ as in \eqref{eq:sec5:example:Control-Policy} is a valid DTCBF for the discrete-time system \eqref{eq:sec5:example:verification:discrete-dynamical}. To this end, we employ and compare the standard $\alpha$BB algorithm \cite{Floudas2010a} with our proposed methods for the verification problem.

Using the standard $\alpha$BB algorithm with convergence tolerance \mbox{$\epsilon_c \defeq 10^{-6}$} and feasibility tolerance $\epsilon_f \defeq 10^{-12}$ (both related to the stopping criteria of the $\alpha$BB algorithm, as detailed in \cite{Floudas2010a}), a global minimizer $x^*$ of the optimization problem \eqref{eq:sec3:known-proposition} is obtained after 55 iterations, each taking an average of $\unit[0.4]{s}$ (total time $\unit[22]{s}$), as 
\begin{align}
    x^*_1 = 0.841, \quad x^*_2 = -1.457. \nonumber
\end{align} 
It is evident that while the DTCBF constraint \eqref{eq:sec2:DTCBF-constraint-alpha} is not satisfied at $x^*$, where $x^* \defeq [x^*_1 ~ x^*_2]^{\T}$, the state $x^*$ lies slightly outside the zero-superlevel set $\mathcal{C}$ of $h$, in the sense that $h(x^*) < 0$. Thus, drawing a conclusion on the validity of $(h,\gamma)$ with the control policy $\pi$ remains inconclusive using the $\alpha$BB algorithm (see Remark \ref{remark:differences-ourmethod-abb}).

Then, we apply Algorithm \ref{algorithm:verification} for the known control policy case with \mbox{$\epsilon_f \defeq \epsilon_h \defeq 10^{-6}$}. After 2 iterations, with each iteration taking an average of $\unit[0.2]{s}$, the output of our method presents a counterexample given by
\begin{align}
\bar{x}_1 = 1.030, \quad \bar{x}_2 = -1.110. \nonumber
\end{align}
It can be observed that the DTCBF constraint \eqref{eq:sec2:DTCBF-constraint-alpha} is not satisfied at $\bar{x}$, where $\bar{x} \defeq [\bar{x}_1 ~ \bar{x}_2]^{\T}$, and  $h(\bar{x}) \geqslant 0$. Thus, either $(h,\gamma)$ is not a valid DTCBF for the system \eqref{eq:sec5:example:verification:discrete-dynamical}, or $\pi$ is not a friend of $(h,\gamma)$.

To determine whether $(h,\gamma)$ is a valid DTCBF with a different control policy, we apply Algorithm \ref{algorithm:verification} for the unknown control policy case with \mbox{$\epsilon_f \defeq \epsilon_h \defeq \epsilon_d \defeq 10^{-6}$}. The application of the proposed algorithm results in Figure \ref{fig:sec5:verification-example}. As observed in Figure~\ref{fig:example-DTCBF:1a}, the subdomains are either outside $\mathcal{C}$ or $(h,\gamma)$ is verified within them. As a result, $(h,\gamma)$ with the computed piecewise constant friend, as depicted in \mbox{Figures \ref{fig:example-DTCBF:1b}} and \ref{fig:example-DTCBF:1c}, is verified as a valid DTCBF after approximately 2500 iterations, each taking an average of $\unit[0.4]{s}$ (total time about $\unit[17]{min}$).

\subsubsection{Cart-Pole} 
\begin{figure}[t]\centering 
\includegraphics[width=0.42\textwidth]{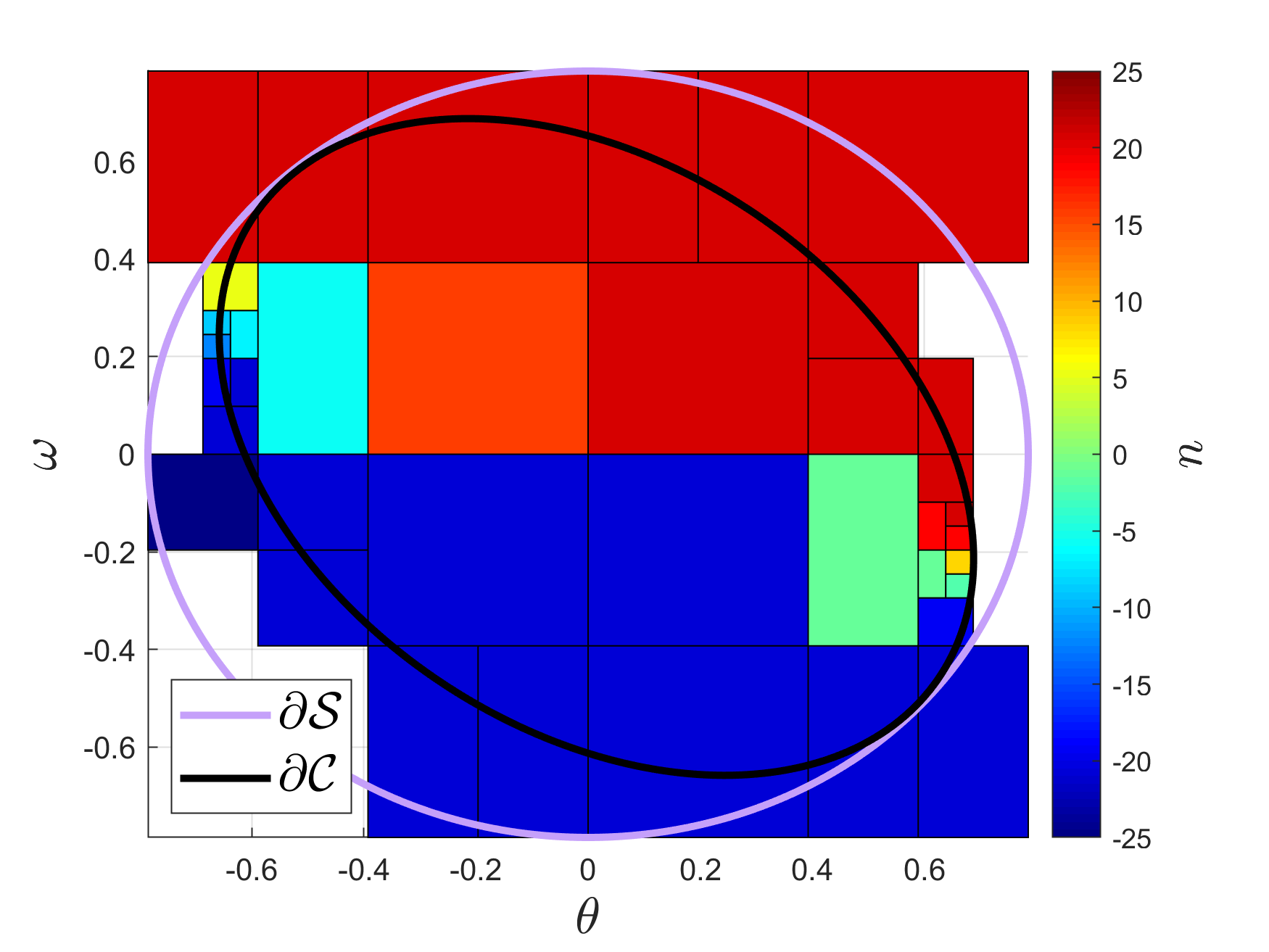}
    \caption{As a cross-check, Algorithm~\ref{algorithm:verification} is applied to the synthesized DTCBF~\eqref{sec5-2:synthesized-DTCBF} with $\gamma(r) \defeq r$, $r \in \mathbb{R}_{\geqslant 0}$, for the discretized cart-pole system~\eqref{eq:dynamics-cart-pole}. The pale purple circle represents the boundary of the safe set $\mathcal{S}$ defined in~\eqref{eq:safe-set-cart-pole}, and the black ellipse depicts the boundary of the zero-superlevel set of the DTCBF. The colored subdomains represent the computed control input values $u$ that satisfy the DTCBF constraint~\eqref{eq:sec2:DTCBF-constraint-alpha} on each corresponding subdomain.}
    \label{fig:sec4:verif-cart-pole-example}
\end{figure}

We again consider the discretized cart-pole system \eqref{eq:dynamics-cart-pole} with the same parameters as mentioned in Section~\ref{sec:sec5:synthesis-cart-pole}. As a cross-check, we aim to verify whether the synthesized $h$ in \eqref{sec5-2:synthesized-DTCBF} with $\gamma(r) = r$, $r \in \mathbb{R}_{\geqslant 0}$, is a valid DTCBF for the system \eqref{eq:dynamics-cart-pole} with the control admissible set~$\mathbb{U}$. First, we apply Algorithm~\ref{algorithm:verification} in the case of a known control policy, using the synthesized control policy $\pi$ given in \eqref{sec5-2:synthesized-control-policy}. After around 160 iterations (approximately $\unit[5]{s}$), the pair $(h, \gamma)$, with its friend $\pi$, is verified as a valid DTCBF. Then, we apply Algorithm~\ref{algorithm:verification} in the case of an unknown control policy, and thus $\pi$ is not assumed to be given. For all subdomains intersecting the zero-superlevel set of~$h$, a control input that satisfies the DTCBF constraint~\eqref{eq:sec2:DTCBF-constraint-alpha} is found in approximately 90 iterations (around~$\unit[8]{s}$). Thus, $(h,\gamma)$ is verified as a DTCBF with a synthesized piecewise constant friend, which differs from the control policy $\pi$ in \eqref{sec5-2:synthesized-control-policy}. This result is illustrated in Figure~\ref{fig:sec4:verif-cart-pole-example}. The values of $\alpha_{\mathcal{F}, i}^{(k)}$ and $\alpha_{\mathcal{H}, i}^{(k)}$, used to construct convex underestimators in the optimization problems \eqref{eq:sec3:verification:known-convex} and \eqref{eq:sec3:verification:unknown:step3-convex}, are computed with the \textit{PyInterval} library.

\section{Conclusion and Future Work}
\label{sec:sec6}
In this paper, we have proposed a novel branch-and-bound method, inspired by the $\alpha$BB algorithm, to either verify a candidate DTCBF as a valid DTCBF for a general nonlinear discrete-time system with input constraints, or falsify it by providing a counterexample (within predefined tolerances). This method is applicable in both cases, whether a corresponding control policy is known or unknown. Additionally, we have presented a novel bilevel optimization approach to synthesize a DTCBF and a corresponding control policy for a general nonlinear discrete-time system with input constraints and an arbitrary safe set. To enhance the tractability of this approach, we have proposed several strategies to lighten the computational burden. We have also applied the proposed methods to numerical case studies.

For future work, we aim to introduce methods tailored to high-dimensional systems for verifying and synthesizing DTCBFs. Moreover, we aim to extend these methods to verify and synthesize \textit{robust} DTCBFs, thereby guaranteeing robust safety for discrete-time systems with disturbances. Our goal is also to adapt the proposed methods for synthesizing discrete-time Control Lyapunov Functions and stabilizing control policies for nonlinear discrete-time systems.

\section*{References}
\bibliographystyle{IEEEtran}
\bibliography{IEEEabrv}

\begin{IEEEbiography}[{\includegraphics[width=1in,height=1.25in,clip,keepaspectratio]{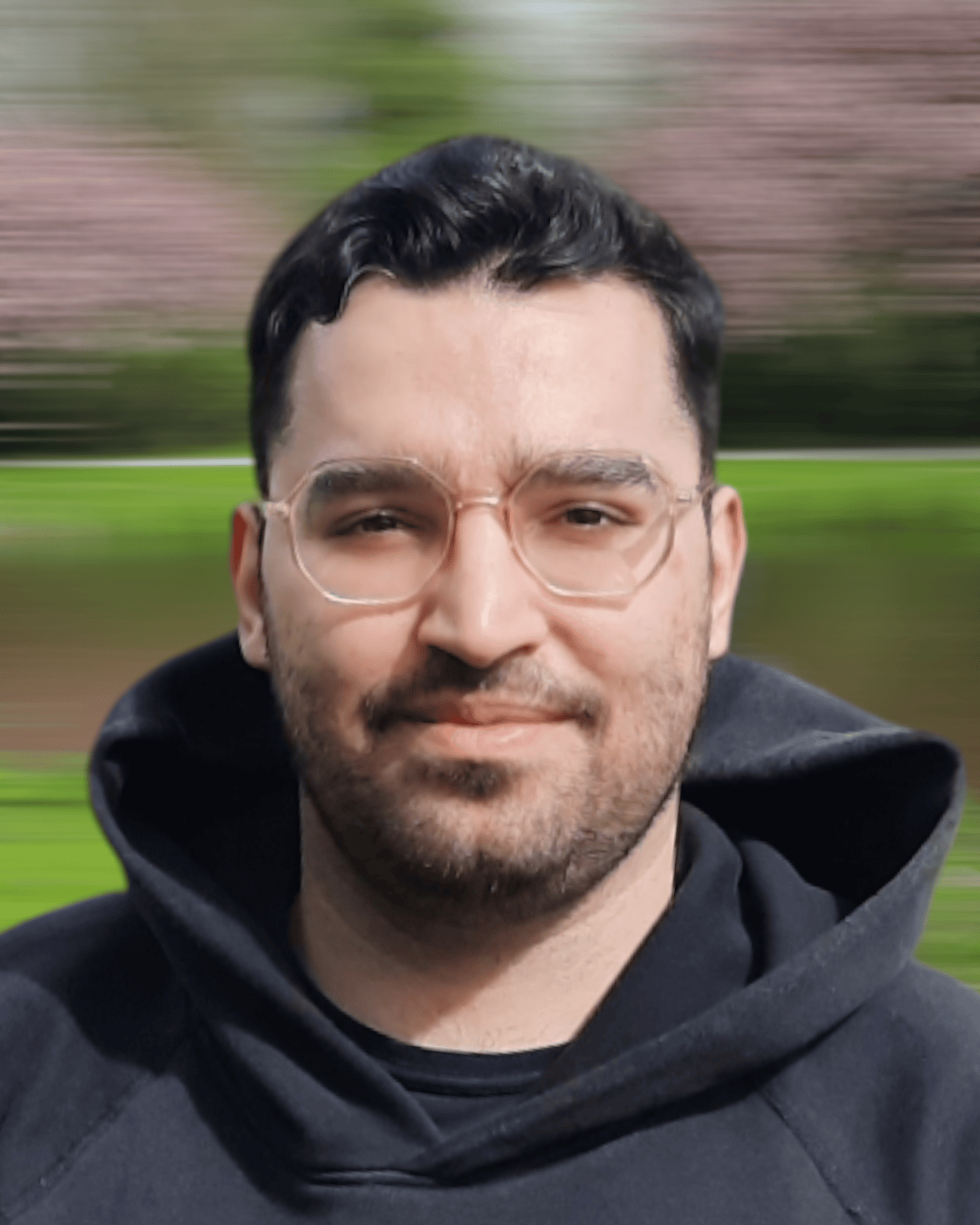}}]{Erfan Shakhesi}(Student Member, IEEE) received the M.Sc. degree (cum laude) in systems and control from the Eindhoven University of Technology (TU/e) in 2023. He is currently pursuing the Ph.D. degree at TU/e. His research interests include nonlinear control, safety-critical control, and model predictive control.
\end{IEEEbiography}

\begin{IEEEbiography}[{\includegraphics[width=1in,height=1.25in,clip,keepaspectratio]{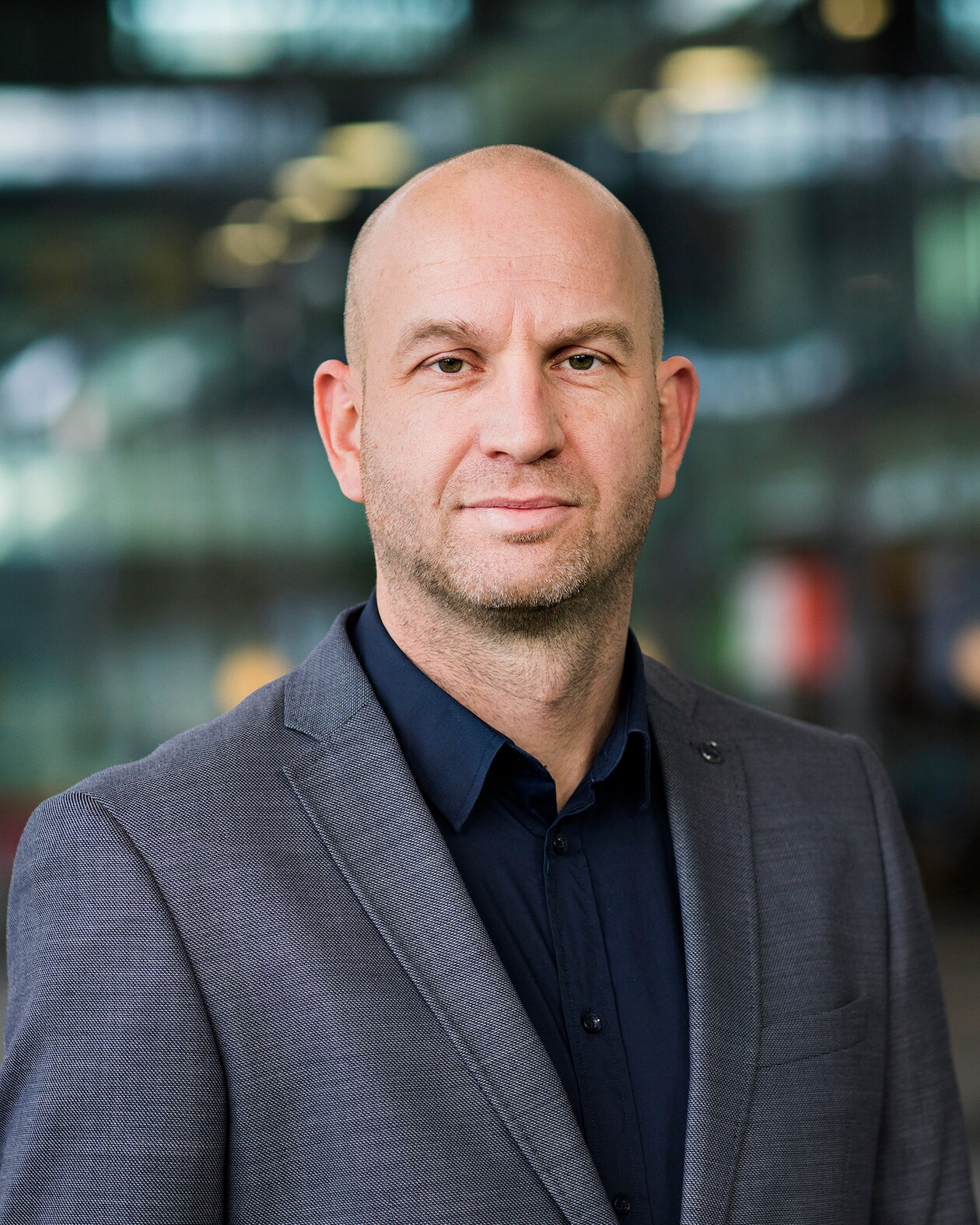}}]{W. P. M. H. (Maurice) Heemels}(Fellow, IEEE) received M.Sc. (mathematics) and Ph.D. (EE, control theory) degrees (summa cum laude) from the Eindhoven University of Technology (TU/e) in 1995 and 1999, respectively. From 2000 to 2004, he was with the Electrical Engineering Department, TU/e, as an assistant professor, and from 2004 to 2006 with the Embedded Systems Institute (ESI) as a Research Fellow. Since 2006, he has been with the Department of Mechanical Engineering, TU/e, where he is currently a Full Professor and Vice-Dean. He held visiting professor positions at ETH, Switzerland (2001), UCSB, USA (2008) and University of Lorraine, France (2020). He is a Fellow of the IEEE and IFAC, and the chair of the IFAC Technical Committee on Networked Systems (2017-2023). He served/s on the editorial boards of Automatica,  Nonlinear Analysis: Hybrid Systems (NAHS), Annual Reviews in Control, and IEEE Transactions on Automatic Control, and is the Editor-in-Chief of NAHS as of 2023. He was a recipient of a personal VICI grant awarded by NWO (Dutch Research Council) and recently obtained an ERC Advanced Grant. He was the recipient of the 2019 IEEE L-CSS Outstanding Paper Award and the Automatica Paper Prize 2020-2022. He was elected for the IEEE-CSS Board of Governors (2021-2023).   His current research includes hybrid and cyber-physical systems, networked and event-triggered control systems and model predictive control. \end{IEEEbiography}

\begin{IEEEbiography}[{\includegraphics[width=1in,height=1.25in,clip,keepaspectratio]{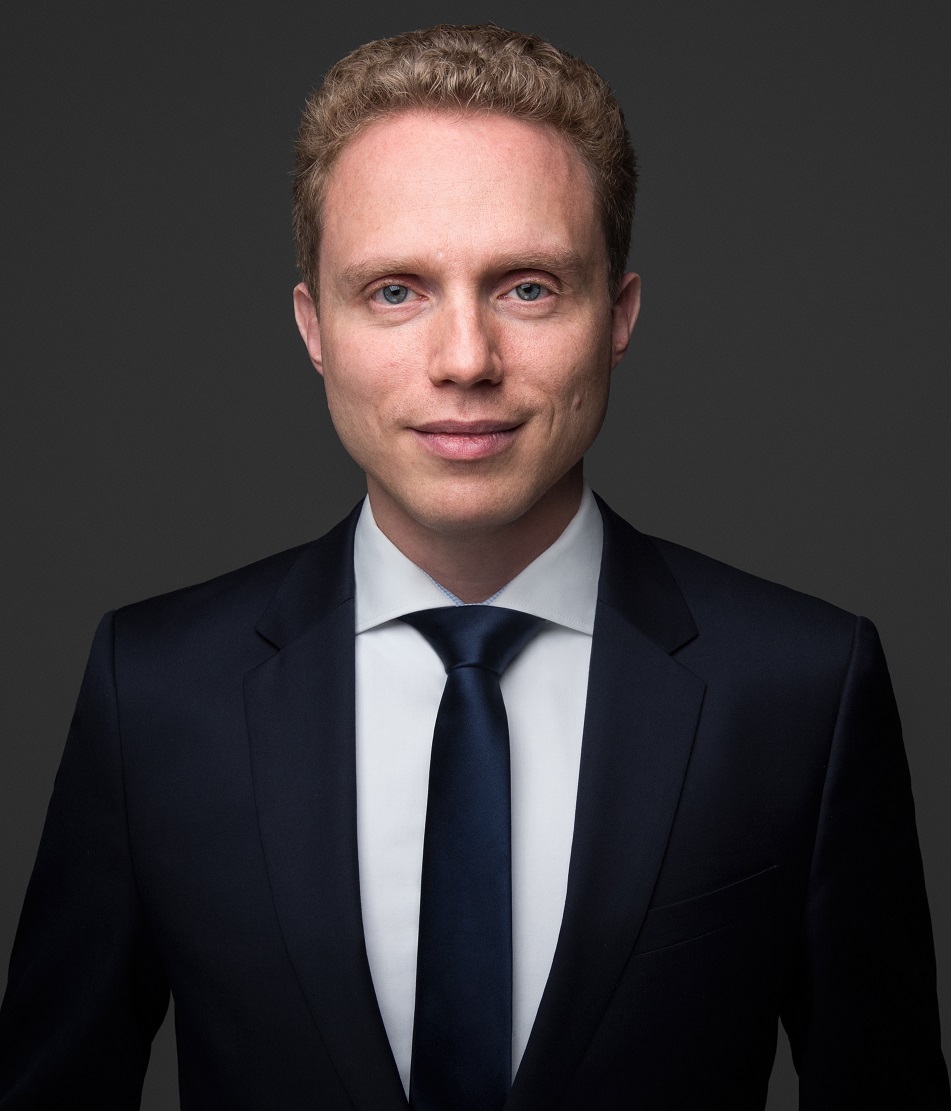}}] {Alexander Katriniok} (Senior Member, IEEE) is an Assistant Professor in the Control Systems
Technology section at the Eindhoven University of Technology (TU/e) in the Netherlands. He
received the Ph.D. degree in Mechanical Engineering from RWTH Aachen University, Germany,
in 2013. His research focuses on learning-based control, provable safety and stability guarantees and real-time numerical optimization, with application to (networked) autonomous systems. He is a Senior Member of IEEE, and the chair of the IEEE Control Systems Society (CSS) Technical Committee on Automotive Controls. He serves as an Associate Editor of the IEEE Transactions on Control Systems Technology (TCST), the IEEE CSS Technology Conference Editorial Board (TCEB), and the IEEE Robotics and Automation Letters (RA-L).
\end{IEEEbiography}

\end{document}